\documentclass{amsart}
\usepackage{amssymb}
\usepackage{latexsym}

\newcommand{\Q}{{\mathbb Q}}
\newcommand{\R}{{\mathbb R}}

\newcommand{\N}{{\mathbb N}}

\newcommand{\Z}{{\mathbb Z}}

\newcommand{\CalW}{{\mathcal{W}}}
\newcommand{\CalA}{{\mathcal{A}}}
\newcommand{\CalB}{{\mathcal{B}}}

\newcommand{\Oomega}{(\Omega,T)}
\newcommand{\OomegaS}{(\Omega (S), T)}

\newcommand{\SL} {\mbox{{\rm SL(2,$\R$)}}  }

\newcommand{\Zet}{ (\mathcal{Z})  }
\newcommand{\SingCont}{ (\mathcal{S}\mathcal{C}) }
\newcommand{\AnTr}{ (\mathcal{A}\mathcal{T}) }

\newtheorem{theorem}{Theorem}
\newtheorem{lemma}{Lemma}[section]
\newtheorem{prop}[lemma]{Proposition}
\newtheorem{coro}{Corollary}
\newtheorem{definition}[lemma]{Definition}
\newtheorem{remark}{Remark}

\begin{document}

\title[Uniform Multiplicative Ergodic Theorem]{A Condition of
Boshernitzan and Uniform Convergence in the Multiplicative Ergodic
Theorem}

\author[D.~Damanik, D.~Lenz]{David Damanik$\,^{1}$, Daniel
Lenz$\,^{2}$}

\thanks{D.\ D.\ was supported in part by NSF grant DMS--0227289}

\maketitle

\vspace{0.3cm}

\noindent $^1$ Department of Mathematics 253--37, California
Institute of Technology, Pasadena, CA 91125, U.S.A., E-Mail:
damanik@its.caltech.edu
\\[0.2cm]
$^2$ Fakult\"at f\"ur Mathematik, TU
Chemnitz, D-09107 Chemnitz, Germany, E-Mail:
dlenz@mathematik.tu-chemnitz.de\\[0.3cm] 2000 AMS Subject
Classification: 37A30, 47B39

\begin{abstract}
This paper is concerned with uniform convergence in the
multiplicative ergodic theorem on aperiodic subshifts. If such a
subshift satisfies a  certain condition,  originally introduced by
Boshernitzan, every locally constant $\SL$-valued cocycle is
uniform. As a consequence, the corresponding Schr\"odinger
operators exhibit Cantor spectrum of Lebesgue measure zero.

An investigation  of Boshernitzan's condition then  shows that  these results cover all
earlier results of this type and, moreover, provide various new
ones. In particular, Boshernitzan's condition is shown to hold for
almost all circle maps and almost all Arnoux-Rauzy subshifts.

\end{abstract}

\section{Introduction}  \label{Introduction}
This paper is concerned with uniform convergence in the multiplicative
ergodic theorem.

More precisely, let $\Oomega$ be a topological dynamical system
Thus, $\Omega$ is a compact metric space and $T : \Omega
\longrightarrow \Omega$ is a homeomorphism. Assume furthermore
that $\Oomega$ is uniquely ergodic, that is, there exists a unique
$T$-invariant probability measure $\mu$ on $\Omega$.

As usual the dynamical system $\Oomega$ is called minimal if every
orbit $\{T^n \omega : n\in \Z\}$ is dense in $\Omega$. It is called
aperiodic if $T^n \omega\neq \omega$ for all $\omega\in \Omega$ and
$n\neq 0$.

Let $\SL$ be the group of real valued $2 \times 2$-matrices with
determinant equal to one equipped with the topology induced by the
standard metric on $2\times 2$ matrices.

To a continuous function $A: \Omega\longrightarrow \SL$ we associate
the cocycle
$$A(\cdot,\cdot):  \Z   \times\Omega \longrightarrow \SL$$
defined by
\begin{equation*}
A(n,\omega) \equiv \left\{\begin{array}{r@{\quad:\quad}l} A(T^{n-1}
 \omega)\cdots A(\omega) & n>0\\ Id & n=0\\ A^{-1} (T^n \omega) \cdots
 A^{-1} (T^{-1}\omega) & n < 0.
\end{array}\right.
\end{equation*}
By the multiplicative ergodic theorem, there exists a $\Lambda(A)\in
\R$ with
\begin{equation}\label{MultET} \Lambda (A)= \lim_{n \to \infty} \frac{1}{n} \log \|
  A(n,\omega)\|
\end{equation}
for $\mu$-almost every $\omega \in \Omega$.  Now, it is well known
that unique ergodicity of $\Oomega$ is equivalent to uniform
convergence in the Birkhoff additive ergodic theorem when applied to
continuous functions. Therefore, it is natural to investigate
uniform convergence in \eqref{MultET}. This motivates the
following definition.

\begin{definition} \cite{Fur, Wal}.   Let $\Oomega$ be uniquely  ergodic.  The continuous function  $A: \Omega \longrightarrow \SL$ is called uniform   if the limit  $\Lambda (A)= \lim_{n \to \infty} \frac{1}{n} \log \| A(n,\omega)\|$ exists  for all $\omega \in \Omega$ and the convergence is uniform on $\Omega$.
\end{definition}

\begin{remark}{\rm  For minimal topological dynamical  systems, uniform existence of the limit in the definition implies uniform convergence. This was   proven by Furstenberg and
    Weiss \cite{FW}. In fact, their result is even more general and
    applies to arbitrary real-valued continuous cocycles.  }
\end{remark}

Various aspects of uniformity of cocycles have been considered in the
past:

\smallskip

A first topic has been to provide examples of non-uniform
cocycles. In fact, in \cite{Wal} Walters asks the question whether
every uniquely ergodic dynamical system  with non-atomic measure
$\mu$      admits  a non-uniform cocycle.  He presents a class of
examples admitting non-uniform cocycles based on results of Veech
\cite{Vee}. He also gives  another  class of examples, namely
suitable irrational rotations, for which  non-uniformity was shown
by Herman \cite{Her}.  In general, however, Walters' question is
still open.

\smallskip

A different line of study has been pursued by Furman in \cite{Fur}. He
 characterizes uniformity of $A$ on a given uniquely ergodic minimal
 $\Oomega$ by a suitable hyperbolicity condition.  The results of
 Furman can essentially be extended to uniquely ergodic systems
 (and, in fact, a strengthening of some sort can be obtained for minimal uniquely
 ergodic systems), as shown by Lenz in \cite{Len4}. They
 also give that the corresponding results of \cite{Her2} provide
 examples of non-uniform cocycles as discussed in \cite{Len4}.

\smallskip

Finally, somewhat complementary to Walters' original question, it
is possible to study conditions on subshifts over finite alphabets
which imply uniformity of locally constant cocycles.  This topic
and variants of it have been discussed at various places
\cite{DL2,Hof,Len1,Len2,Len3, Len4}. It is the main focus of the
present article.  It is not only of intrinsic interest but also
relevant in the study of spectral theory of certain Schr\"odinger
operators, as recently shown by Lenz \cite{Len2} (see below for
details).

\smallskip

To elaborate on this and state our main results, we recall some further
notions.

$\Oomega$ is called a subshift over $\CalA$ if $\CalA$ is finite
with discrete topology and $\Omega$ is a closed $T$-invariant
subset of $\CalA^\Z$, where $\CalA^\Z$ carries the product
topology and $T : \CalA^\Z \longrightarrow \CalA^\Z$ is given by
$(T s ) (n) := s (n+1)$.  A function $F$ on $\Omega$ is called
locally constant if there exists an $N\in  \N$ with
\begin{equation} \label{LocallyConstant}
F(\omega) = F(\rho)\;\:\mbox{whenever}\;\:   (\omega(-N),\ldots,
\omega(N))=(\rho(-N),\ldots, \rho(N)).
\end{equation}

We will freely use notions from combinatorics on words (see, e.g.,
\cite{Lot1, Lot2}). In particular, the elements of $\CalA$ are
called letters and the elements of the free monoid $\CalA^\ast$
over $\CalA$ are called words.  The length $|w|$ of a word $w$ is
the number of its letters.  The number of occurrences of a word
$w$ in a word $x$ is denoted by $\#_w (x)$.

Each subshift $\Oomega$ over $\CalA$ gives rise to the associated set of
words
\begin{equation}
\label{subwordsomega}
\CalW (\Omega) :=\{ \omega(k) \cdots \omega(k + n -1) : k\in \Z, n\in
\N, \omega \in \Omega\}.
\end{equation}

For $w \in \CalW$, we define
\begin{equation*}
V_w :=\{\omega \in \Omega : \omega (1) \cdots \omega (|w|) = w\}.
\end{equation*}

Finally, if $\nu$ is a $T$-invariant probability measure on $\Oomega$ and $n\in \N$, we set
\begin{equation}
\eta_\nu (n) :=\min  \{\nu (V_w) : w\in \CalW, |w|=n\}.\end{equation}
If $\Oomega$ is uniquely ergodic with invariant probability measure $\mu$ , we set $\eta (n):= \eta_\mu (n)$.

\begin{definition} Let $\Oomega$ be a subshift over $\CalA$. Then,
  $\Oomega$ is said to satisfy condition {\rm (B)} if there exists  an ergodic
    probability  measure $\nu$ on $\Omega$ with
$$\limsup_{n\to \infty} n \,  \eta_\nu (n) >0.$$
Thus, $\Oomega$ satisfies {\rm (B)} if and only if  there exists
an ergodic probability  measure $\nu$ on $\Omega$,  a constant
$C>0$ and a sequence $(l_n)$ in $\N$ with $l_n \to \infty$ for
$n\to \infty$ such that $|w| \nu( V_w)  \geq C $ whenever $w\in
\CalW (\Omega)$ with $|w|= l_n$ for some $n\in \N$.
\end{definition}

This condition was introduced by Boshernitzan in \cite{Bosh1}
(also see \cite{Bosh4} for related material). For minimal interval
exchange transformations, it was shown to imply unique ergodicity
by Veech in \cite{Vee2}. Finally, in \cite{Bosh2}, Boshernitzan
showed that it implies unique ergodicity for arbitrary minimal
subshifts.

\smallskip

Our main result is:

\begin{theorem}\label{main} Let $\Oomega$ be a minimal subshift which satisfies {\rm (B)}.
Let $A : \Omega \longrightarrow \SL$ be locally constant. Then, $A$ is uniform.
\end{theorem}

As discussed below, this result covers all  earlier results of this
form as given in \cite{DL2, Hof,Len3,Len4}. Moreover, as we will show
below, it also applies to various new examples, including many circle
maps and Arnoux-Rauzy subshifts. This point is worth emphasizing, as most circle maps and Arnoux-Rauzy subshifts seem to have been rather out of reach of earlier methods.

The proof of the main result is based on two steps. In the first
step, we give various equivalent characterizations of condition
(B). This is made precise in Theorem \ref{Characterization}. This
result may be of independent interest. In our context it shows
that (B) implies uniform convergence on ``many scales.''  In the
second step, we use the so-called Avalanche Principle introduced
by Goldstein and Schlag in \cite{GS} and extended by Bourgain and
Jitomirskaya in  \cite{BJ}  to conclude uniform convergence from
uniform convergence on ``many scales.''

As a by-product of our proof, we obtain a simple combinatorial
argument for unique ergodicity for subshifts satisfying (B).
Unlike the proof given in \cite{Bosh2}, we do not need any apriori
estimates on the number of invariant measures.

\medskip

As mentioned already, our results are particularly relevant in the
study of certain Schr\"odinger operators. This is discussed next:

To each bounded  $V: \Z \longrightarrow\R$, we can associate the Schr\"odinger
operator $H_V : \ell^2 (\Z)\longrightarrow \ell^2 (\Z)$  acting by
$$(H_V u )(n) \equiv u(n+1) + u(n-1) + V(n) u(n).$$
The spectrum of $H_V$ is denoted by $\sigma(H_V)$.

Now, let  $\Oomega$ be a subshift over $\CalA$ and assume without
loss of generality that $\CalA \subset \R$. Then, $\Oomega$ gives
rise to the  family $(H_\omega)_{\omega\in \Omega}$ of selfadjoint
operators. These operators arise in the study of aperiodically
ordered solids, so-called quasicrystals. They exhibit interesting
spectral features such as Cantor spectrum of Lebesgue measure
zero, purely singularly continuous spectrum and anomalous
transport.  They have attracted a lot of attention in recent years
(see, e.g., the surveys \cite{Dam, Sut3} and discussion below for
details). Recently, Lenz has shown that uniformity of certain
locally constant cocycles is intimately related to Cantor spectrum
of Lebesgue measure zero for these operators \cite{Len2}. This can
be combined with our main result to give the following theorem
(see below for details).

\begin{theorem}\label{zeromeasure} Let $\Oomega$ be a minimal subshift
  which satisfies {\rm (B)}.  If $\Oomega$ is aperiodic, then there
  exists a Cantor set $\Sigma\subset \R$ of Lebesgue measure zero with
  $\sigma(H_\omega) = \Sigma$ for every $\omega \in \Omega$.
\end{theorem}

This result covers  all earlier results on Cantor spectrum of
measure zero \cite{AD,BBG,BIST,BG,DL4,DL5,Len2,LTWW,OL2,Sut,Sut2} as
discussed below. More importantly, it gives various new ones. In
particular, it covers almost all circle maps and Arnoux-Rauzy
subshifts.

\smallskip

To give a flavor of these new examples, we mention the following
theorem. Define for $\alpha, \theta, \beta \in (0,1)$ arbitrary,
the function
$$V_{\alpha,\beta,\theta} : \Z\longrightarrow \{0,1\},\;\;\mbox{
 by}\;\: V_{\alpha,\beta,\theta} (n):=\chi_{ [1-\beta,1)} (n\alpha +
 \theta \mod 1),$$ where $\chi_M$ denotes the characteristic function
 of the set $M$. These functions are called circle maps.

\begin{theorem}\label{CantorCircleMap}
Let $\alpha \in (0,1)$ be irrational. Then, we have the following:
\\[1mm]
{\rm (a)} For almost every $\beta\in (0,1)$, the spectrum
$\sigma(H_{V_{\alpha,\beta,\theta }})$ is a Cantor set of Lebesgue
measure zero for every $\theta \in (0,1)$.
\\[1mm]
{\rm (b)} If $\alpha$ has bounded continued fraction expansion,
then $\sigma(H_{V_{\alpha,\beta,\theta }})$ is a Cantor set of
Lebesgue measure zero for every $\beta\in (0,1)$ and every $\theta
\in (0,1)$.
\end{theorem}

\begin{remark}{\rm
This result is particularly relevant as all earlier results on
Cantor spectrum for circle maps \cite{AD,BIST,DL4,Sut,Sut2} only
cover a set of parameters $(\alpha,\beta)$ of Lebesgue measure
zero in $(0,1)\times (0,1)$ (cf.~Appendix~\ref{highpowers}). }
\end{remark}

Finally, we mention the following by-product of our investigation.
Details (and precise definitions) will be discussed in Section
\ref{Application}.

\begin{theorem} \label{ContinuityLE} Let $\Oomega$ be a minimal subshift
  which satisfies {\rm (B)} and $(H_\omega)_{\omega\in \Omega}$ the associated family of operators. Then the Lyapunov exponent $\gamma : \R \longrightarrow [0,\infty)$ is continuous.
\end{theorem}

\medskip

The paper is organized as follows: In Section \ref{Key} we study
condition (B) and show its equivalence to various other
conditions. As a by-product this shows unique ergodicity of
subshifts satisfying (B). Moreover, it is used in
Section~\ref{Uniformity} to give  a proof of our main result.
Stability of the results under certain operations on the subshift
is discussed in  Section~\ref{Stability}. In Section \ref{Known}
we discuss examples for which (B) is known to hold. New examples,
viz certain circle maps and Arnoux-Rauzy subshifts, are given in
Sections~\ref{Circle} and \ref{Arnoux-Rauzy}. Finally, the
application to Schr\"odinger operators is discussed in Section
~\ref{Application}.

\section{Boshernitzan's Condition (B)} \label{Key}
In this section, we give various equivalent characterizations of
(B). This is made precise in Theorem \ref{Characterization}. Then, we
provide a new proof of unique ergodicity for systems satisfying (B) in
Theorem \ref{UE}. Theorem \ref{Characterization} in some sense
generalizes the main results of \cite{Len1} and its proof heavily uses
and extends ideas from there.

\medskip

To state our result, we need some preparation. We start by introducing
a variant of Boshernitzan's condition (B).  Namely, if $\Oomega$ is a
subshift, we define for $w\in \CalW (\Omega)$ the set $U_w$ by
$$U_w :=\{ \omega\in \Omega : \exists n\in \{0,1,\ldots, |w|-1\}
\text{ such that } \omega (-n + 1)\ldots  \omega (-n + |w|) =
w\}.$$ If $\omega$ belongs to $U_w$, we say that $w$ occurs in
$\omega$ around one.

\begin{definition} Let $\Oomega$ be a subshift over $\CalA$. Then,
$\Oomega$ is said to satisfy condition {\rm (B')} if there exists
an ergodic probability  measure $\nu$ on $\Omega$, a constant
$C'>0$, and a sequence $(l'_n)$ in $\N$ with $l'_n \to \infty$ for
$n\to \infty$ such that $\nu( U_w) \geq C'$ whenever $w\in \CalW
(\Omega)$ with $|w|= l'_n$ for some $n\in \N$.
\end{definition}

Next, we discuss a consequence of Kingman's ergodic theorem.
Recall that $F : \CalW (\Omega)\longrightarrow \R$ is called
subadditive if it satisfies $F(x y) \leq F(x) + F(y)$ whenever
$x,y,xy\in \CalW (\Omega)$, where $\Oomega$ is an arbitrary
subshift.

\begin{prop}\label{defLambda} Let $\Oomega$ be a uniquely ergodic subshift with invariant
probability measure $\mu$.  Let $F : \CalW (\Omega)\longrightarrow \R$
be subadditive, then there exists a number $\Lambda(F)\in \R\cup \{-\infty\}$ with
$$\Lambda(F) = \lim_{n\to \infty} n^{-1} F(\omega(1)\cdots \omega(n))$$
for $\mu$-almost every
$\omega$ in $\Omega$.
\end{prop}
\begin{proof} For $n\in \N$, we define the continuous function $f_n : \Omega\longrightarrow \R$, by
$$ f_n (\omega) := F(\omega(1)\ldots \omega(n)).$$ As $F$ is subadditive, $(f_n)$  is  a subadditive cocycle.
Thus Kingman's subadditive theorem applies. This  proves the
statement.
\end{proof}

\begin{theorem} \label{Characterization} Let $\Oomega$ be a minimal  subshift over $\CalA$. Then the
  following conditions are equivalent:
\begin{itemize}
\item[(i)]  $\Oomega$ satisfies {\rm (B)}. \item[(ii)] $\Oomega$
satisfies {\rm (B')}. \item[(iii)]   $\Oomega$ is uniquely ergodic
and there exists  a sequence $(l'_n)$ in $\N$ with $l'_n \to
\infty$ for $n\to \infty$ such that $\lim_{n\to \infty} |w_n|^{-1}
F(w_n) = \Lambda (F) $  for every subadditive $F$ and every
sequence $(w_n)$ in $\CalW (\Omega)$ with $|w_n|= l'_n$ for every
$n\in \N$.
\end{itemize}
\end{theorem}

The remainder of this section is devoted to a proof of this
theorem. The proof will be split into several parts.

\begin{lemma}
\label{BvsBprime} Let $\Oomega$ be a minimal subshift. Then,
$\Oomega$ satisfies {\rm (B)} if and only if it satisfies {\rm
(B')}.
\end{lemma}
\begin{proof}
If $\Oomega$ is periodic, validity of (B) and (B') is immediate. Thus,
we can restrict our attention to aperiodic $\Oomega$.

\smallskip

Apparently, $\nu (U_w) \leq |w| \nu (V_w)$ for all $w\in \CalW (\Omega)$ and
all ergodic probability measures $\nu$ on $\Omega$. Thus, (B') implies
(B) (with the same $\nu$, $l_n$, and $C$).

\smallskip

Conversely, assume that $\Oomega$ satisfies (B). We will show that
it satisfies (B') with $l'_n = [2 l_n / 3 ] + 1$, $n\in \N$,  and
$C'= C/9$. Here, for arbitrary $a\in \R$, we set  $[a] := \sup\{
n\in Z : n\leq a\}$.

Consider $v\in \CalW (\Omega)$ with $|v| = l'_n$ for some $n\in \N$. Choose
$w\in \CalW (\Omega)$ with $|w|= l_n$ such that $v$ is a prefix of $w$. There
are two cases:

{\it Case 1.} There exists a primitive $x\in \CalW (\Omega)$ and a prefix
$\widetilde{x} $ of $x$ such that $ w = x^k \widetilde{x}$ for some
$k\geq 6$.

\smallskip

As $\Oomega$ is minimal and aperiodic, the word $x$ does not occur
with arbitrarily high powers. Thus, we can find $y\in \CalW
(\Omega)$ such that
$$ \widetilde{w} := x^{k-1} y\in \CalW (\Omega)$$ satisfies $ |\widetilde{w}|=
l_n$ but $x^k$ is not a prefix of $\widetilde{w}$. Now, as $x$ is
primitive, it does not appear non-trivially in $x^{k-1}$.
Therefore, different copies of $\widetilde{w}$ have distance at
least $(k-2) |x|$. This gives
$$ \nu (U_{ \widetilde{w}} ) \geq (k-2 ) |x| \nu (V_{\widetilde{w} } )
\geq \frac{ (k-2) |x|}{ (k +2) |x|} | \widetilde{w}|  \nu
(V_{\widetilde{w} } ) \geq \frac{1}{2} C.$$ Moreover, by construction,
$v$ is a subword of $\widetilde{w}$ (and even of $x^{k-1}$) with
$$ \frac{|v|}{|\widetilde{w}|} \geq \frac{1}{2}.$$
Putting these estimates together, we infer
$$ \nu (U_v) \geq \frac{1}{2} \nu (U_{\widetilde{w} }) \geq
\frac{1}{2} \cdot \frac{1}{2} \cdot C = \frac{C}{4}.$$

{\it Case 2.} There does not exist a primitive $x$ in $\CalW$ and a
prefix $\widetilde{x} $ of $x$ with $w = x^k \widetilde{x}$ for some
$k\geq 6$.

\smallskip

In this case, different copies of $w$ have distance at least
$\frac{1}{6} |w|$. Therefore, we have
$$ \nu (U_w) \geq \frac{1}{6} |w| \nu (V_w)$$
and this gives
$$ \nu (U_v) \geq \frac{2}{3} \nu (U_w) \geq \frac{2}{3}\cdot \frac{1}{6}\cdot |w| \nu (V_w)\geq \frac{1}{9} C.$$

\smallskip

In both cases the desired estimates hold and the proof of the lemma is
finished.
\end{proof}

We next give our proof of unique ergodicity for systems satisfying
(B'). The proof proceeds in two steps. In the first step, we use
(B') to show existence of the frequencies along certain sequences.
In the second step, we show existence of the frequencies along all
sequences. Let us emphasize that it is exactly this two-step
procedure which is underlying the proof of our main result on
locally constant matrices. However, in that case the details are
more involved.

\smallskip

We need the following proposition.

\begin{prop}\label{Birkhoff} Let $\Oomega$ be a subshift with ergodic probability  measure $\nu$. Let $f :\Omega \longrightarrow \R$ be a bounded measurable function. Then,
$$\lim_{n,m\geq 0, n+ m \to \infty} \frac{1}{n + m} \sum_{k=-m}^n
f(T^k \omega) = \nu (f)$$ for $\nu$-almost every $\omega\in \Omega$.
\end{prop}
\begin{proof}  By Birkhoff's ergodic theorem, we have  both
$$\lim_{n\to \infty} \frac{1}{n} \sum_{k=0}^{n-1} f(T^k \omega) = \nu
(f)\;\: \mbox{and}\;\: \lim_{m\to \infty} \frac{1}{m} \sum_{k=-m}^{0}
f(T^k \omega) = \nu(f)$$ for $\nu$-almost every $\omega\in
\Omega$. Now, for every sequence $(a_k)_{k\in \Z}$ with
$$\lim_{n\to \infty} \frac{1}{n}\sum_{k= 0}^n a_k = \lim_{m\to \infty}
 \frac{1}{m}\sum_{k=-m }^0 a_k = a,$$ one easily infers
 $$\lim_{n,m\geq 0, n + m\to \infty} \frac{1}{n+m} \sum_{k=-m}^n a_k =
 a.$$ The statement follows immediately.
\end{proof}

\begin{theorem} \label{UE} If the subshift $\Oomega$ satisfies {\rm (B')},
it is uniquely ergodic and minimal.
\end{theorem}
\begin{proof} It suffices to show that the frequencies $\lim_{|x|\to \infty}
\frac{ \#_w (x)}{|x|}$ exist for every $w\in \CalW$. Then, the system is
uniquely ergodic by standard reasoning. Moreover, in this case, the system
is minimal as well as all frequencies are positive by (B').

\smallskip

Thus, let an arbitrary $w\in \CalW (\Omega)$ be given. We proceed in two steps.

{\it Step 1.} For all $\varepsilon>0$, there exists an $n_0 = n_0
(\varepsilon)$ with $\left| \frac{\#_w (x)}{|x|} - \nu (V_w)\right|
\leq \varepsilon$ whenever $|x|=l'_n$ with $n\geq n_0$.

\smallskip

{\it Step 2.} For $\varepsilon>0$, there exists an $N_0 = N_0
(\varepsilon)$ with $\left| \frac{\#_w (x)}{|x|} - \nu (V_w)\right|
\leq \varepsilon$ whenever $|x|\geq N_0$.

\smallskip

Here, Step 2 follows easily from Step 1 by partitioning long words $x$
into pieces of length $l'_{n}$ with sufficiently large $n\in \N$.

Thus, we are left with the task of proving Step 1. To do so, assume
the contrary. Then, there exist $\delta >0$, $(x_n)$ in $\CalW$ and
$(l'_{k(n)})$ in $\N$ with $ |x_n|= l'_{k(n)}$, $k(n) \longrightarrow
\infty$ and
\begin{equation} \label{groesser} \left| \frac{\#_w (x_n)}{|x_n|} - \nu (V_w)\right| \geq \delta
\end{equation}
for every $n\in \N$.  Consider
$$E:= \bigcap_{n=1}^\infty \bigcup_{k=n}^\infty U_{x_k}.$$
By (B'), we have
$$ \nu (E) = \lim_{n\to \infty} \nu (\cup_{k=n}^\infty U_{x_k} ) \geq C'>0.$$
Thus, by Proposition \ref{Birkhoff}, we can find an $\omega $ in $E$ with
\begin{equation} \label{kleiner}
\lim_{n,m\geq 0, n + m\to \infty} \frac{\#_{w} (\omega (-m) \ldots
\omega (n))}{n + m} = \nu (V_w).
\end{equation}
As $\omega$ belongs to $E$, there are infinitely many $x_n$
occurring around one  in $\omega$. Now, if we calculate the
occurrences of $w$ along this sequence of $x_n$, we stay away from
$\nu (V_w)$ by at least $\delta$ according to \eqref{groesser}. On
the other hand, by \eqref{kleiner}, we come arbitrarily close to
$\nu (V_w)$ when calculating the frequency of $w$ along any
sequence of words occurring in $\omega$ around one. This
contradiction proves Step~1 and therefore finishes the proof of
the theorem by the discussion above.
\end{proof}

Our next task is to relate (B') and convergence in subadditive ergodic
theorems. We need two auxiliary results.

\begin{prop}\label{obereSchranke} Let $\Oomega$ be a uniquely ergodic
subshift and $F : \CalW (\Omega) \longrightarrow \R$ be subadditive.
Then, $\limsup_{|x|\to \infty} |x|^{-1} F(x) \leq \Lambda (F)$.
\end{prop}
\begin{proof} Define $f_n$ as in the proof of Proposition \ref{defLambda}. Then,
the statement is a direct consequence of Corollary~2 in
\cite{Fur}.
\end{proof}

\begin{prop}\label{glm}
Let $\Oomega$ be a uniquely ergodic subshift with invariant
probability measure $\mu$. Let $w\in \CalW(\Omega)$ be arbitrary and
denote by $\chi_{U_w}$ the characteristic function of $U_w$.  Then,
$$\lim_{n\to\infty} \frac{1}{n} \sum_{k=0}^{n-1} \chi_{U_w} (T^k
\omega) = \mu (U_w)$$ uniformly in $\omega\in  \Omega$.
\end{prop}
\begin{proof}
As $U_w$ is both closed and open, the characteristic function
$\chi_{U_w}$ is continuous. Thus, the statement follows from unique
ergodicity.
\end{proof}

Now, our result on subadditive ergodic theorems and (B') reads as
follows.

\begin{lemma}\label{SET}
Let $\Oomega$ be a uniquely ergodic and minimal subshift. Let $(w_n)$
be a sequence in $\CalW (\Omega)$ with $|w_n|\longrightarrow \infty$,
$n\to \infty$. Then, the following assertions are equivalent:

\begin{itemize}

\item[(i)] $\lim_{n\to \infty} |w_n|^{-1} F(w_n) = \Lambda (F)$ for
every subadditive $F : \CalW (\Omega) \longrightarrow \R$.

\item[(ii)] There exists a $C'>0$ with $\mu (U_{w_n}) \geq C'$ for every
$n\in \N$.
\end{itemize}
\end{lemma}
\begin{proof} The proof can be thought of as an adaptation and extension
of the proofs of Lemma~3.1 and Lemma~3.2 in \cite{Len1} to our
setting.

\smallskip

(i) $\Longrightarrow$ (ii). Assume the contrary. Then, the sequence
$(\mu (U_{w_n}))$ is not bounded away from zero. By passing to a
subsequence, we may then assume without loss of generality that
\begin{equation}\label{bb}
\sum_{n=1}^\infty \mu (U_{w_n}) < \frac{1}{2}.
\end{equation}
As $\Oomega$ is minimal, we have $\mu (U_{w_n}) >0$ for every $n\in
\N$. Moreover, by assumption, we have

\begin{equation}\label{bbb}
|w_n|\longrightarrow \infty, n\longrightarrow \infty.
\end{equation}
For $w,x\in\CalW (\Omega)$, we say that $w$ occurs in $x$ around
$j\in\{1,\ldots,|x|\}$ if there exists $l\in \N$ with $l \leq j <
l + |w| -1$ and $ x (l) \ldots x (l + |w| -1) = w$.

Now, define for $n \in \N$, the function $F_n : \CalW
(\Omega) \longrightarrow \R$ by
$$ F_n (x) := \#\{ j\in \{1,\ldots, |x|\} : \mbox{$w_n$ occurs in
$x$ around $j$}\}.$$  Here, $\# M$ denotes the cardinality of $M$. Thus $F_n (x)$ measures the amount of ``space''
covered in $x$ by copies of $w_n$.  Obviously, $- F_n$ is  subadditive for every $n\in \N$.

\smallskip
The definition of $F_n$ shows
$$ F_n (\omega(1)\ldots \omega (m) ) = \sum_{k=0}^{m- |w_n| - 1} \chi_{U_{w_n}}
(T^k \omega)$$ for arbitrary $\omega\in \Omega$ and $m\in \N$. Thus,
by Proposition \ref{glm}, we have
$$\lim_{|x|\to \infty} |x|^{-1} F_n (x) = \mu (U_{w_n})$$
for arbitrary but fixed $n\in \N$.

Invoking this equality and \eqref{bb} and \eqref{bbb}, we can choose
inductively for every $k\in \N$ a number $n(k)\in \N$ with

$$\frac{| w_{n(k+1)}|  }{2} > |w_{n(k)}|$$
and
$$ \sum_{j=1}^k \frac{F_{n(j)} (x)}{|x| } < \frac{1}{2},$$ whenever
$|x|\geq |w_{n(k+1)}|$. It is not hard to see
that
$$ F(x) :=\sum_{j=1}^\infty F_{n(2 j)} (x)$$
 is finite for every $x\in \CalW (\Omega)$  and  $-F :\CalW (\Omega) \longrightarrow \R$, $x\mapsto - F(x)$,   is
subadditive. Therefore, by our assumption (i) the limit
$$- \Lambda(- F) = \lim_{n\to \infty} \frac{F(w_n)}{|w_n|}$$
exists. On the other
hand, for every $k\in \N$, we have
$$\frac{F (w_{n(2k)})}{|w_{n(2k)}|} \geq \frac{F_{n(2k)}
(w_{n(2k)})}{|w_{n(2k)}|} =1$$ as well as
$$\frac{F (w_{n(2k+1)})}{|w_{n(2k+1)}|} = \frac{1}{| w_{n(2k+1)}|}
\sum_{j=1}^k F_{n(2j)} (w_{n(2k+1)})\leq \frac{1}{| w_{n(2k+1)}|}
\sum_{j=1}^{2k}  F_{n(j)} (w_{n(2k+1)}) < \frac{1}{2}.$$ This is a
contradiction and the proof of this part of the lemma  is finished.

\medskip

(ii) $\Longrightarrow$ (i). Let $F : \CalW (\Omega)\longrightarrow \R$
be subadditive. By Proposition \ref{obereSchranke}, we have
\begin{equation}\label{aa}
\limsup_{|x|\to \infty} \frac{F(x)}{|x|} \leq \Lambda (F).
\end{equation}
Thus, it remains to show
$$\Lambda(F) \leq \liminf_{n\to \infty} \frac{F(w_n)}{|w_n|}.$$ Assume
the contrary.  Then, $\Lambda(F) > -\infty$ and there exists a
subsequence $(w_{n(k)})$ of $(w_n)$ and $\delta >0$ with
\begin{equation} \label{uu} \frac{F(w_{n(k)})}{|w_{n(k)}|} \leq \Lambda(F) - \delta
\end{equation}
for every $k\in \N$. For $w,x\in \CalW (\Omega)$, we define
$\#_w^\ast (x)$ to be the maximal number of disjoint copies of $w$
in $x$.

\smallskip

It is not hard to see that
$$|w|\cdot \#_w^\ast (\omega (1)\ldots \omega(m)) \geq
\frac{1}{2} \sum_{k=0}^{m -|w| -1} \chi_{U_w} (T^k \omega)$$ for all
$\omega\in \Omega$ and $m\in \N$.  By Proposition \ref{glm}, this implies
$$\liminf_{|x|\to \infty} \frac{\#_w^\ast (x) }{|x|} |w|\geq
\frac{1}{2} \mu ( U_w).$$

\smallskip
Combining this with our assumption (ii), we infer
\begin{equation} \label{stern}
\liminf_{|x|\to \infty} \frac{\#_{w_{n(k)}}^\ast (x)}{|x|}
|w_{n(k)}| \geq \frac{C'}{2}
\end{equation}
for every $k\in \N$. By \eqref{aa}, we can choose $L_0$ such that
\begin{equation} \label{oben}
\frac{F(x)}{|x|} \leq \Lambda (F) + \frac{C'}{16} \delta,
\end{equation}
whenever $|x|\geq L_0$.  Fix $k\in \N$ with $|w_{n(k)}|\geq L_0$.
Using \eqref{stern}, we can now find an $ L_1\in\R$ such that
every $x\in \CalW (\Omega)$ with $|x|\geq L_1$ can be written as $x= x_1
w_{n(k)} x_2 w_{n(k)} \ldots x_l w_{n(k)} x_{l+1}$ with
\begin{equation}\label{hatata}
\frac{l-2}{2} \geq \frac{C'}{8} \frac{|x|}{|w_{n(k)}|}.
\end{equation}
Now, considering only every other copy of $w_{n(k)}$ in $x$, we can
write $x$ as $ x=y_1 w_{n(k)} y_2 \ldots y_r w_{n(k)} y_{r+1}$, with
$|y_j|\geq |w_{n(k)}|\geq L_0$, $j=1,\ldots,r+1,$ and by
\eqref{hatata}
$$r\geq \frac{l-2}{2} \geq \frac{C'}{8} \frac{ |x|}{ |w_{n(k)}|}.$$
Using $\eqref{oben}$, $\eqref{uu}$ and this estimate, we can now
calculate
\begin{eqnarray*}
\frac{F(x)}{|x|} &\leq & \sum_{j=1}^{r+1} \frac{ F(y_j)}{|y_j|}
\frac{|y_j|}{|x|} + \frac{F(w_{n(k)})}{|w_{n(k)}|}\frac{ r
|w_{n(k)}|}{|x|}\\ &\leq & \sum_{j=1}^{r+1} (\Lambda (F)
+\frac{C'}{16}\delta )\frac{|y_j|}{|x|} + ( \Lambda (F) -\delta) \frac{
r |w_{n(k)}|}{|x|}\\ &\leq& \Lambda (F) + \frac{C'}{16}\delta -
\frac{C'}{8} \frac{ |x|}{ |w_{n(k)}|}\frac{|w_{n(k)}|}{|x|} \delta \\
&\leq & \Lambda (F) - \frac{C'}{16}\delta.
\end{eqnarray*}

As this holds for arbitrary $x\in \CalW (\Omega)$ with $|x|\geq L_1$, we arrive at
the obvious contradiction $\Lambda (F) \leq \Lambda (F) - \frac{C'}{16}\delta
$. This finishes the proof.
\end{proof}

\begin{proof}[Proof of Theorem \ref{Characterization}.]
Given the previous results, the proof is simple: The equivalence
of (i) and (ii) is shown in Lemma \ref{BvsBprime}. The implication
(ii) $\Longrightarrow$ (iii) follows from Theorem \ref{UE}
combined with Lemma \ref{SET}. The implication (iii)
$\Longrightarrow$ (ii) is immediate from Lemma \ref{SET}. This
finishes the proof of Theorem \ref{Characterization}.
\end{proof}

\section{Uniformity of Locally Constant Cocycles}\label{Uniformity}
In this section we provide a proof of our main result, Theorem
\ref{main}. As mentioned already, the cornerstones of the proof
are Theorem \ref{Characterization} and the so-called Avalanche
Principle, introduced in \cite{GS} and later extended in
\cite{BJ}.

\medskip

We use the Avalanche Principle in the following form given in
Lemma 5 of \cite{BJ}.
\begin{lemma} \label{Avalanche}
There exist constants $\lambda_0 >0$ and $\kappa >0$ such that
$$ \left| \log\| A_N \ldots A_1\| + \sum_{j=2}^{N-1} \log \|A_j\| -
\sum_{j=1}^{N-1} \log \|A_{j+1 } A_{j}\| \right| \leq \frac{\kappa \cdot N
}{\exp (\lambda)},$$ whenever $N= 3^P$ with $P\in \N$ and $A_1,\ldots,
A_N$ are elements of $\SL$ such that

\begin{itemize}

\item $\log\|A_j\| \geq \lambda \geq \lambda_0$ for every $j=1,\ldots, N$;

\item $| \log \|A_j\| + \log\|A_{j+1}\|- \log\|A_j A_{j+1}\| | <\frac{1}{2} \lambda$  for every $j=1,\ldots, N$.

\end{itemize}
\end{lemma}
\begin{remark}{\rm Actually,  Lemma 5 in \cite{BJ} is  more general in that more general $N$ are allowed.}\end{remark}

Before we can  give the proof of Theorem \ref{main}, we need one more   auxiliary result.

\begin{prop}\label{Vergleich} Let $\Oomega$ be an arbitrary subshift and $A : \Omega \longrightarrow \SL$  a locally constant function. Then,
$$ 0 = \lim_{n\to \infty} \sup \left\{ \frac{1}{n} \left| \log \|
A(n,\omega)\| - \log\| A(n,\rho)\|\right| : \omega (1)\ldots
\omega (n) = \rho(1)\ldots \rho(n) \right\}.$$
\end{prop}

\begin{proof} As $A$ is locally constant, there exists an $N\in \N$ such that $A(\omega) = A (\rho)$, whenever $\omega (-N)\ldots \omega (N) = \rho(-N) \ldots \ \rho (N)$. Thus,
$$
A(n- 2N,T^N \omega) = A (n - 2N, T^N \rho),
$$ whenever $n\geq 2N$ and $\omega (1)\ldots \omega(n) = \rho(1)
\ldots \rho (n)$.  Moreover, for arbitrary matrices $X,Y,Z$ in
$\SL$, we have
$$ \log\|Y\| - \log \|X\| - \log\|Z\| \leq \log\| X Y Z\| \leq
\log\|X\| + \log \|Y\| + \log \|Z\|,
$$ where we used the triangle inequality as well as $\|M\| = \|
M^{-1}\|$ for $M\in \SL$.  Finally, we have
$$
 A (n,\sigma) = A (N, T^{n-N} \sigma) A (n - 2N, T^N\sigma) A(N, \sigma).$$
Putting these three equations together, we arrive at the desired conclusion.
\end{proof}

\begin{remark}{\rm Let us point out that the previous proposition is the only point in our considerations where local constancy of $A$ enters. In particular, our main result holds for all $A$ for which the conclusion of the proposition holds.  }\end{remark}

\begin{proof}[Proof of Theorem \ref{main}.]
Let $\Oomega$ be a subshift satisfying (B) and let $A : \Omega
\longrightarrow \SL$ be locally constant. We have to show that $A$
is uniform.

\smallskip

\textit{Case 1.} $\Lambda(A) =0$: As $A$ takes values in $\SL$, we
have $\|A(n,\omega)\|\geq 1$ and the estimate
$$0 \leq \liminf_{n \to \infty} \frac{1}{n} \log\|A(n,\omega)\|$$
holds uniformly in $\omega\in \Omega$. On the other hand, by Corollary
2 of \cite{Fur}, we have
$$\limsup_{n\to \infty} \frac{1}{n} \log\|A(n,\omega)\| \leq
\Lambda(A) $$ uniformly in $\omega\in \Omega$. This shows the desired
uniformity in this case.

\medskip

\textit{Case 2.} $\Lambda(A) >0$: Define $F : \CalW (\Omega)
\longrightarrow \R$ by
$$ F(x) := \sup\{ \log\|A (n,\omega)\| : \omega (1)\ldots \omega(n) =
x\}.$$ Apparently, $F$ is subadditive. As discussed above, there
exists then $\Lambda(F)$ with
$$\Lambda(F) = \lim_{n\to \infty} \frac{F(\omega(1)\ldots \omega
(n))}{n}$$ for $\mu$-almost every $\omega\in \Omega$. On the other
hand, by the multiplicative ergodic theorem, there also exists
$\Lambda(A)$ with
$$ \Lambda (A) = \lim_{n\to \infty} \frac{\log\|A(n,\omega)\|}{n}$$
for $\mu$-almost every $\omega\in \Omega$. By Proposition \ref{Vergleich},
we infer that $ \Lambda (A) = \Lambda (F).$ Summarizing, we have
\begin{equation} \label{AgleichF}
\Lambda(A) = \Lambda(F) >0.
\end{equation}
Combining this equation with Theorem \ref{Characterization}, we infer
$$ \lim_{n\to \infty} \frac{F(w_n)}{|w_n|} =\Lambda (A),$$
whenever $(w_n)$ is a sequence with $|w_n| = l'_n$.
Also, combining \eqref{AgleichF} with Proposition \ref{obereSchranke}, we infer
$$ \limsup_{n\to \infty} \frac{1}{n} \log\|A(n,\omega)\| \leq
\limsup_{|x|\to \infty} \frac{F(x)}{|x|} \leq \Lambda (A)$$ uniformly
in $\omega \in \Omega$. It remains to show
$$\Lambda(A) \leq \liminf_{n\to \infty} \frac{1}{n} \log\|A(n,\omega)\|$$
uniformly in $\omega\in \Omega$. To do so, let  $\varepsilon >0$ with $\varepsilon \leq 1/12$ be given.

The preceding considerations and Proposition \ref{Vergleich} give
existence of $n_0\in \N$ such that with
$$l := \frac{l_{n_0}'}{2},$$ the following holds:

\begin{itemize}
\item[(I)] $\log\| A(n,\omega)\| \leq \Lambda (A) (1+ \varepsilon) n$ for
all $\omega\in \Omega$ whenever $n\geq l$.

\smallskip

\item[(II)] $\log\|A(2 l,\omega)\| \geq \Lambda(A) (1 - \varepsilon) 2l$
for all $\omega\in \Omega$.

\smallskip

\item[(III)] $ \Lambda(A) (1 - 3\varepsilon) l \geq \lambda_0$.

\smallskip

\item[(IV)]  $\frac{2 \kappa}{ l\exp (\lambda_0)}  < \varepsilon \Lambda (A)$.

\end{itemize}
Here, $\lambda_0$ and $\kappa$ are the constants from Lemma
\ref{Avalanche}. Using (II), subadditivity and (I), we can
calculate
\begin{eqnarray*}
\Lambda (A) (1 -\varepsilon) 2 l &\leq & \log\|A( 2l, \omega)\|\\
&\leq &\log\|A(l,\omega)\| + \log\|A (l, T^l \omega)\|\\
&\leq&  \log\|A(l,\omega)\| + \Lambda(A) (1 + \varepsilon) l.
\end{eqnarray*}
This implies $\Lambda(A) (1 - 3 \varepsilon) l \leq \log\|A(l,\omega)\|$
and therefore by (III), \begin{equation}\label{Voraus1} \lambda_0 \leq
\Lambda(A) (1 - 3 \varepsilon) l \leq \log\|A(l,\omega)\|
\end{equation}
for every $\omega\in \Omega$.
Moreover, by subadditivity, (I)  and (II), we have
\begin{eqnarray*}
& &\left| \log\| A(l,\omega)\| + \log\|A(l,T^l\omega)\| - \log\|A(2
l,\omega)\|\right|\\ &=& \log\| A(l,\omega)\| + \log\|A(l,T^l\omega)\|
- \log\|A(2 l,\omega)\| \\ &\leq& \Lambda(A) 2 l (1 + \varepsilon) -
\log\|A(2 l,\omega)\| \\ &\leq & \Lambda(A) 2 l (1 + \varepsilon) -
\Lambda(A) 2 l (1 -\varepsilon)\\ &=& \Lambda(A) 4 l \varepsilon
\end{eqnarray*}
for arbitrary $\omega \in \Omega$. Using the assumption $\varepsilon \leq 1/12$, we infer
\begin{equation}\label{Voraus2}
\left| \log\| A(l,\omega)\| + \log\|A(l,T^l\omega)\| -
\log\|A(l,\omega)\|\right| \leq \frac{1}{2} \Lambda(A) (1 - 3
\varepsilon) l.
\end{equation}
Equations \eqref{Voraus1} and \eqref{Voraus2} and (III) show that
the Avalanche Principle, Lemma \ref{Avalanche},  with
$$\lambda =\Lambda(A) (1 -3  \varepsilon) l$$
 can be applied to the matrices $A_1,\ldots, A_N$,
where $N= 3^P$ with $P\in \N$ arbitrary  and
$$A_j = A(l, T^{(j-1)l} \omega),\;\:j=1,\ldots, N$$
with $\omega\in \Omega$ arbitrary. This gives
$$ \left| \log\| A_N \ldots A_1\| + \sum_{j=2}^{N-1} \log \|A_j\| -
\sum_{j=1}^{N-1} \log \|A_{j+1 } A_{j}\| \right| \leq \frac{\kappa N
}{\exp (\lambda)}.$$
This yields

\begin{eqnarray*}
\log\|A_N\ldots A_1\| &\geq& \sum_{j=1}^{N-1} \log\|A_{j+1} A_j\| -
\sum_{j=2}^{N-1} \log\|A_j\| - \frac{\kappa \cdot N }{\exp
(\lambda)}\\ &\geq & (N-1) \Lambda(A) (1 -\varepsilon) 2 l - (N-2)
\Lambda(A) (1+ \varepsilon) l - \frac{\kappa \cdot N }{\exp (\lambda)}\\
&= & \Lambda(A) N l (1 - 3\varepsilon) + \Lambda(A) 4 \varepsilon l -
\frac{\kappa \cdot N }{\exp (\lambda)}\\ &\geq & \Lambda(A) N l (1 -
3\varepsilon) - \frac{\kappa \cdot N }{\exp (\lambda)}.
\end{eqnarray*}
Here, we used (I) and (II) in the second step and positivity of
$\Lambda(A) 4 \varepsilon l $ in the last step. Dividing by by $n:= N l$,
and invoking (IV), we obtain
\begin{equation}\label{fastfertig}
\Lambda(A) (1- 4 \varepsilon) \leq \frac{1}{n} \log\|A(n,\omega)\|
\end{equation}
 for all $\omega\in\Omega$ and all $n =  3^P\cdot l $ with  $P\in
 \N$.

\smallskip

We finish the proof by showing that

\begin{equation}\label{fertig}
\Lambda(A) (1 - 44 \varepsilon) \leq \frac{1}{n} \log\|A(n,\omega)\|
\end{equation}
for all $n\geq l $ and all $\omega\in \Omega$. As $\varepsilon$
was arbitrary, this gives the desired statement.

To show  \eqref{fertig}, choose $\omega\in \Omega$ and  and $n\geq l $. Let $P\in \N\cup\{0\}$ be such that
$$ 3^P\cdot l  \leq n <  3^{P+1}\cdot l .$$
Then, by \eqref{fastfertig} and subadditivity we have
\begin{eqnarray*}
 \Lambda(A) (1 - 4\varepsilon)& \leq &  \frac{1 }{  3^{P+2} l  }\log\| A( 3^{P+2} l ,
\omega)\| \\
&\leq & \frac{1 }{  3^{P+2} l  } \log \|A(n,\omega)\| + \frac{1 }{  3^{P+2} l  }\ \log\| A( 3^{P +2} l  - n,
T^n\omega)\|\\
&\leq& \frac{1}{n} \log\|A(n,\omega)\| \cdot  \frac{n}{  3^{P+2} l } + \Lambda(A)  (1 + \varepsilon) ( 1 - \frac{ n}{ 3^{P+2} l }),
\end{eqnarray*}
where we could use (I) in the last estimate as, by assumption on $n$, $ 3^{P +2} l  - n \geq 3^{P+1} 2 l  >l$.
Now, a direct calculation gives
$$\Lambda(A) \left(1 + \varepsilon - 5\varepsilon \frac{ 3^{P+2}l}{n}\right) \leq \frac{1}{n} \log\|A(n,\omega)\|.$$
As   $3^{P+2} l / n \leq 9$ by the very choice of $P$,  the
desired equation \eqref{fertig} follows easily. This finishes  the
proof of our main theorem.
\end{proof}

\section{Stability of Uniform Convergence Under Substitutions} \label{Stability}
In the last section, we studied sufficient conditions on $\Oomega$
to ensure property
$$ {\rm (P)}:\;\:\mbox{Every locally constant $A : \Omega \longrightarrow
\SL$ is uniform.} $$ In this section, we consider ``perturbations''
$\OomegaS$ of $\Oomega$ by substitutions $S$ and study how validity of
(P) for $\Oomega$ is related to validity of (P) for $\OomegaS$.

\medskip

We start with the necessary notation. Let $\CalA$ and $\CalB$ be
finite sets. A map $S :\CalA \longrightarrow \CalB^\ast$ is called
a substitution.  Obviously, $S$ can be extended to $\CalA^\ast$ in
the obvious way. Moreover, for a two-sided infinite word
$(\omega(n))_{n\in \Z}$ over $\CalA$, we can define $S(\omega)$ by
$$S(\omega) := \cdots S(\omega(-2)) S(\omega(-1)) | S(\omega(0))
S(\omega(1)) S(\omega(2))\cdots,$$ where $|$ denotes the position of
zero. If $\Oomega$ is a subshift over $A$ and $S : \CalA \longrightarrow
\CalB^\ast$ is a substitution, we define $\Omega(S)$ by
$$\Omega(S) : =\{T^k S(\omega) : \omega \in \Omega, k\in \Z\}.$$ Then,
$\OomegaS$ is a subshift over $\CalB$. It is not hard to see that
$\OomegaS$ is minimal (uniquely ergodic) if $\Omega$ is minimal
(uniquely ergodic).

\begin{theorem}\label{StabilSubst} Let $\Oomega$ be a minimal uniquely ergodic
subshift over $\CalA$ that satisfies {\rm (P)}. Let $S$ be a
substitution over $\CalA$. Then, $\OomegaS$ satisfies {\rm (P)} as
well.
\end{theorem}

\begin{proof} Let $ B :\Omega (S) \longrightarrow \SL$ be locally constant. Define
$$ A :\Omega\longrightarrow \SL\;\:\mbox{by}\;\: A (\omega) := B(
|S(\omega(0))|,S(\omega)).$$ Then, $A$ is locally constant as well and
$$A(n,\omega) = B(|S (\omega(0)\ldots \omega(n-1))  |, S(\omega)).$$
In particular, we have
\begin{equation}\label{crux}
\frac{\log\| B(|S (\omega(0)\ldots \omega(n)) |, S(\omega))\| }{|
S(\omega(0)\ldots \omega(n))|} = \frac{n + 1}{ | S(\omega(0)\ldots
\omega(n))|} \cdot \frac{\log\|A(n,\omega)\|}{n + 1 }.
\end{equation}
By
$$|S(\omega(0)\ldots \omega (n)) | = \sum_{a\in \CalA} |S(a)|
\#_a (\omega(0)\ldots \omega(n))$$
and unique ergodicity of
$\Oomega$, the quotients
$$ \frac{n + 1 }{ | S(\omega(0)\ldots \omega(n))|}$$ converge uniformly
in $\omega \in \Omega$ towards a number $\rho$. From \eqref{crux} and
validity of (P) for $\Oomega$ we infer that
$$\lim_{n\to \infty} \frac{ \log\|
B(|S (\omega(0)\ldots \omega(n)) |, S(\omega))\| }{| S(\omega(0)\ldots \omega(n))|}  = \rho  \cdot \Lambda(A)$$
uniformly on $\Omega$. As every $\sigma\in \Omega(S)$ has the form
$\sigma = T^k S(\omega)$ with $|k| \leq \max\{ |S(a)| : a\in \CalA\}$,
uniform convergence of $\frac{1}{n} \log \|B(n,\sigma)\|$ follows.
\end{proof}

In certain cases, a converse of this theorem holds. To be more
precise, let $\Oomega$ be a subshift over $\CalA$ and $S$ a
substitution on $\CalA$. Then, $S$ is called \textit{recognizable}
(with respect to $\Oomega$) if there exists a locally constant map
$$\widetilde{S} : \Omega (S) \longrightarrow \Omega\times \Z$$ with $
\widetilde{S}(T^k S(\omega)) = (\omega,k)$, whenever $0\leq k \leq
|S(\omega (0))|$.
Recognizability is known for various classes of substitutions that
generate aperiodic subshifts, including all primitive
substitutions \cite{Mos} and all substitutions of constant length
that are one-to-one \cite{App} (cf.~the discussion in
\cite{Fogg}).

\begin{theorem}\label{StabilRuecksubst} Let $\Oomega$ be a uniquely ergodic
minimal subshift over $\CalA$. Let $S$ be a recognizable
substitution over $\CalA$. If $\OomegaS$ satisfies {\rm (P)}, then
$\Oomega$ satisfies {\rm (P)} as well.
\end{theorem}

\begin{proof} Let $B : \Omega \longrightarrow\SL$ be locally constant. For $\sigma\in\Omega(S)$ define
\begin{equation*}
A(\sigma) \equiv \left\{\begin{array}{r@{\quad:\quad}l} B(\omega) &
 \sigma = S(\omega) \\ id & \mbox{otherwise}.
\end{array}\right.
\end{equation*}
Note that $\sigma = S(\omega)$ if and only if the second component of
$\widetilde{S} (\sigma) $ is $0$. As $\widetilde{S}$ is locally
constant, this shows that $A$ is locally constant as well.

Moreover, by definition of $A$ and recognizability of $S$, we have
$$ A(|S(\omega(0)\ldots \omega(n-1))|, S(\omega)) = B(n,\omega).$$
Now, the proof can be finished similarly to the proof of the previous theorem.
\end{proof}

There is an instance of the previous theorem that deserves special
attention, viz subshifts derived by return words. Return words and
the derived subshifts have been discussed by various authors since
they were first introduced by Durand in \cite{Du}. We recall the
necessary details next.

 Let $\Oomega$ be a minimal subshift and $w\in
\CalW(\Omega)$ arbitrary. Then, $x\in \CalW (\Omega)$ is called a
return word of $w$ if $xw$ satisfies the following three properties:
it belongs to $\CalW(\Omega)$, it starts with $w$ and it contains
exactly two copies of $w$. We then introduce a new alphabet $\CalA_w$
consisting of the return words of $w$. Obviously, there is a natural
map
$$S_w : \CalA_w \longrightarrow \CalA^\ast$$ which maps the return
word $x$ of $w$ (which is a letter of $\CalA_w$) to the word $x$
over $\CalA$. Partitioning every word $\omega\in\Omega$ according
to occurrences of $w$, we obtain a unique two-sided infinite word
$\omega_w$ over $\CalA_w$ with
$$ T^{-k} S_w (\omega_w) = \omega$$
for $k \leq 0$ maximal  with
$\omega(k)\ldots \omega(k + |w|-1) = w$. We define
$$ \Omega_w :=\{\omega_w : \omega \in \Omega\}.$$ Then, $(\Omega_w,T)$
is a subshift, called the subshift derived from $\Oomega$ with
respect to $w$. It is not hard to see that $(\Omega_w,T)$ is
minimal. Moreover, $(\Omega_w,T)$ is uniquely ergodic if $\Oomega$
is uniquely ergodic. Clearly, $S_w$ is recognizable and $\Oomega =
(\Omega_w (S_w),T)$ since the whole construction only depends on
the (local) information of occurrences of $w$. Thus, we obtain the
following corollary from the previous theorem.

\begin{coro} Let $\Oomega$ be a minimal uniquely ergodic subshift that satisfies {\rm (P)}.
Let $w\in \CalW(\Omega)$ be arbitrary. Then, $(\Omega_w,T)$
satisfies {\rm (P)} as well.
\end{coro}

The aim of this paper is to study (P). Given that (B) is a
sufficient condition for (P), it is then natural to ask for
stability properties of (B) as well. It turns out that (B) shares
the stability features of (P).

\begin{theorem}\label{StabilityB} Let $\Oomega$ be a minimal
uniquely ergodic subshift over $\CalA$. Let $S$ be a substitution
on $\CalA$ and $ (\Omega(S), T)$ the corresponding subshift.
\\[1mm]
{\rm (a)} If $\Oomega$ satisfies {\rm (B)}, so does $(\Omega(S),
T)$.
\\[1mm]
{\rm (b)} If $(\Omega(S), T)$ satisfies {\rm (B)} and $S$ is
recognizable, then $\Oomega$ satisfies {\rm (B)} as well.
\end{theorem}

Before we can give a proof, we note the following simple observation.

\begin{prop}\label{Beobachtung}
Let $\Oomega$ be a  minimal uniquely ergodic subshift satisfying
{\rm (B)} with length scales $(l_n)$ and constant $C>0$. Then,
$$ |w| \mu (V_w) \geq \frac{C}{N},$$ whenever $w\in \CalW (\Omega)$
satisfies $l_n / N \leq |w| \leq l_n$ for some $n\in \N$ and $N\in N$.
\end{prop}

\begin{proof} Every $w\in \CalW (\Omega)$  with $l_n / N \leq |w|\leq l_n$
is a  prefix of a $v\in \CalW$ with $|v|= l_n$.  Then, $V_v
\subset V_w$ holds and (B) implies
$$|w|\mu (V_w) \geq \frac{|v|}{N} \mu (V_w) \geq \frac{|v|}{N} \mu
(V_v) \geq \frac{C}{N}.$$ This finishes the proof of the proposition.
\end{proof}

\begin{proof}[Proof of Theorem \ref{StabilityB}.]
Define  $M:=\{|S(a)| : a \in \CalA\}$ and  denote the unique
$T$-invariant probability measure on $\Omega$ (resp., $\Omega(S)$)
by $\mu$ (resp., $\mu_S$).
\\[1mm]
(a) We assume that $\Oomega$ satisfies (B) with length scales
$(l_n)$ and constant $C>0$. Let $w\in \CalW (\Omega(S))$ with $|w|
= l_n$ for some $n \in \N$ be given. Then, there exists a word
$v\in \CalW(\Omega)$ such that $w$ is a subword of $S(v)$ and
satisfies the estimate
\begin{equation}\label{vauzuwe}
 \frac{|w|}{M} \leq |v|\leq |w|.
\end{equation}
Choose $\omega \in\Omega$ arbitrary.   Obviously,
$$\#_w (S(\omega(1)\ldots \omega(k))) \geq \#_v (\omega(1)\ldots \omega(k)).$$
Thus, counting occurrences of $w\in S(\omega)$ and occurrences of $v$ in $\omega$, we obtain by unique ergodicity
\begin{eqnarray*}
|w| \mu_S (V_w)& =& |w|\lim_{n\to \infty} \frac{ \#_w  ( S(\omega)
(1)\ldots S(\omega) (n))}{n}
\geq   |w| \lim_{k\to \infty} \frac{ \#_{v}  (\omega (1) \ldots \omega (k))}{ k M} \\
&=& \frac{|w|}{M} \mu (V_v)
\geq    \frac{1}{M} |v| \mu (V_v)
\geq    \frac{1}{M^2} C,
\end{eqnarray*}
where we used \eqref{vauzuwe} in the second-to-last step and
Proposition \ref{Beobachtung} combined with \eqref{vauzuwe} in the
last step. This shows (B) for $(\Omega(S),T)$ along the same
length scales $(l_n)$ with new constant $C / M^2$.
\\[1mm]
(b) We assume that $(\Omega(S), T)$ satisfies (B) with constant $C
> 0$ and length scales $(l_n)$. By recognizability, there exists a map $\widetilde{S} : \Omega (S)
\longrightarrow \Omega\times \Z$ and an $N\in \N$ with $
\widetilde{S}(T^k S(\omega)) = (\omega,k)$, whenever $0\leq k \leq
|S(\omega (0))|$, and $\widetilde{S} (\omega) = \widetilde{S}
(\rho)$, whenever $\omega(-N) \ldots \omega(N) = \rho (-N)\ldots
\rho (N)$. Let $n_0$ be chosen such that
$$\left[ \frac{l_n}{3 M}\right] \geq N,$$ for all $n\geq n_0$.

Choose an arbitrary $v\in \CalW (\Omega)$ with $|v|= \left[
\frac{l_n}{3 M}\right]$ for some $n\geq n_0$.

Let $x,y \in \CalW (\Omega)$ be given with $|x|= |y| = |v|$ and $x
v y \in \CalW (\Omega)$.  By recognizability and our choice of the
lengths of $x,y$ and $v$, occurrences of $S(x v y)$ in $S(\omega)$
correspond to occurrences of $v$ in $\omega$ for any $\omega\in
\Omega$. Thus, we obtain

$$\#_v (\omega(1)\ldots \omega(n)) \geq \#_{S(x v y)}
(S(\omega(1)\ldots \omega(n)))$$ for every $n\in \N$ and $\omega
\in \Omega$. Therefore, a short calculation invoking unique
ergodicity gives
\begin{eqnarray*}
|v| \mu ( V_v) & =& |v| \lim_{n\to \infty} \frac{ \#_v
  (\omega(1)\ldots \omega(n)) }{n} \geq  |v| \lim_{n\to \infty} \frac {
  \#_{S(x v y)} (S(\omega(1)\cdots \omega(n))) } {n} \\ &\geq &
\frac{|v|}{|S(xvy)|} \lim_{n\to \infty} \frac{|S(\omega(1)\cdots
  \omega (n))| }{n} |S(x v y)| \frac{\#_{S(x v y)}
  (S(\omega(1)\cdots \omega(n))) } {|S(\omega(1)\cdots \omega (n))|}\\
&\geq & \frac{1}{3 M} |S(x v y)| \mu_S (V_{S(xvy)}),
\end{eqnarray*}
where we used the trivial bound $|S(x)| / |x|\geq 1$ in the
second-to-last step. By construction, we have
$$ \frac{l_n}{2 M} \leq |x v y| \leq |S(x v y)| \leq l_n.$$ Thus, we
can apply Proposition \ref{Beobachtung}, and the assumption (B) on
$\Omega(S)$, to our estimate on $|v| \mu(V_v)$ to obtain $|v| \mu
(V_v) \geq \frac{C}{6 M^2 }$. As $v\in \CalW$ with $|v| =\left[
\frac{l_n}{3 M}\right]$ was arbitrary, we infer (B) with the new
length scales $[l_n / 3] $ for $n\geq n_0$ and new constant $ C/
(6 M^2)$.
\end{proof}

\section{Examples Known to Satisfy (B) }\label{Known}

In this section we discuss the classes of subshifts for which the
Boshernitzan condition is either known or a simple consequence of
known results. In our discussion of the occurrence of zero-measure
Cantor spectrum for Schr\"odinger operators in
Section~\ref{Application}, this will be relevant since all the
models for which this spectral property was previously known will
be shown to satisfy (B). Hence we present a unified approach to
all these results.

\subsection{Examples Satisfying (PW): Linearly Recurrent Subshifts and
Subshifts Generated by Primitive Substitutions}

A subshift $\Oomega$ over $\CalA$ satisfies the condition (PW)
(for \textit{positive weights}) if there exists a constant $C > 0$
such that
$$
\liminf_{|x| \to \infty} \frac{\#_v(x)}{|x|} |v| \ge C \text{ for
every } v \in \CalW (\Omega).
$$
This condition was introduced by Lenz in \cite{Len1}. There, it
was shown that the class of subshifts satisfying (PW) is exactly
the class of subshifts for which a uniform subadditive ergodic
theorem holds. Moreover, (PW) implies minimality and unique
ergodicity.

The following is obvious:

\begin{prop}\label{pwimpliesb}
If the subshift $\Oomega$ satisfies {\rm (PW)}, then it satisfies
{\rm (B)}.
\end{prop}

The condition (PW) holds in many cases of interest. For example,
it is easily seen to be satisfied for all linearly recurrent
subshifts. Here, a subshift $\Oomega$ is called \textit{linearly
recurrent} (or \textit{linearly repetitive}) if there exists a
constant $K$ such that if $v,w \in \CalW (\Omega)$ with $|w| \ge K
|v|$, then $v$ is a subword of $w$.

We note:

\begin{prop}\label{lrimpliespw}
If the subshift $\Oomega$ is linearly recurrent, then it satisfies
{\rm (PW)}.
\end{prop}

The class of linearly recurrent subshifts was studied, for
example, in \cite{Du2,DHS}.

A popular way to generate linearly recurrent subshifts is via
primitive substitutions. A substitution $S : \CalA \to \CalA^*$ is
called \textit{primitive} if there exists $k \in \N$ such that
for every $a,b \in \CalA$, $S^k(a)$ contains $b$. Such a
substitution generates a subshift $\Oomega$ as follows. It is easy
to see that there are $m \in \N$ and $a \in \CalA$ such that
$S^m(a)$ begins with $a$. If we iterate $S^m$ on the symbol $a$,
we obtain a one-sided infinite limit, $u$, called a
\textit{substitution sequence}. $\Omega$ then consists of all
two-sided sequences for which all subwords are also subwords of
$u$. One can verify that this construction is in fact independent
of the choice of $u$, and hence $\Omega$ is uniquely determined by
$S$. Prominent examples are given by
$$
\begin{array}{|l|l|}
\hline
a \mapsto ab, \; b \mapsto a & \text{ Fibonacci}\\
\hline
a \mapsto ab, \; b \mapsto ba & \text{ Thue-Morse}\\
\hline
a \mapsto ab, \; b \mapsto aa & \text{ Period Doubling}\\
\hline a \mapsto ab, \; b \mapsto ac, \; c \mapsto db, \; d
\mapsto dc & \text{ Rudin-Shapiro}\\
\hline
\end{array}
$$

The following was shown in \cite{DHS}:

\begin{prop}\label{psimplieslr}
If the subshift $\Oomega$ is generated by a primitive
substitution, then it is linearly recurrent.
\end{prop}

It may happen that a non-primitive substitution generates a
linearly recurrent subshift. An example is given by $a \mapsto
aaba$, $b \mapsto b$. In fact, the class of linearly recurrent
subshifts generated by substitutions was characterized in
\cite{DL5}. In particular, it turns out that a subshift generated
by a substitution is linearly recurrent if and only if it is
minimal.

\subsection{Sturmian and Quasi-Sturmian Subshifts}

Consider a minimal subshift $\Oomega$ over $\CalA$. Recall that
the associated set of words is given by
$$
\CalW (\Omega) :=\{ \omega(k) \cdots \omega(k + n -1) : k\in \Z,
n\in \N, \omega \in \Omega\}.
$$
The (factor) complexity function $p : \N \to \N$ is then
defined by
\begin{equation}\label{complexitydef}
p(n) = \# \CalW_n (\Omega) ,
\end{equation}
where $ \CalW_n (\Omega) = \CalW (\Omega) \cap \CalA^n$ and $\#$ denotes cardinality.

It is a fundamental result of Hedlund and Morse that periodicity
can be characterized in terms of the complexity function
\cite{hm}:

\begin{theorem}[Hedlund-Morse]
$\Oomega$ is aperiodic if and only if $p(n) \ge n+1$ for every $n
\in \N$.
\end{theorem}

Aperiodic subshifts of minimal complexity, $p(n) = n+1$ for every
$n \in \N$, exist and they are called Sturmian. If the complexity
function satisfies $p(n) = n + k$ for $n \ge n_0$, $k,n_0 \in
\N$, the subshift is called quasi-Sturmian. It is known that
quasi-Sturmian subshifts are exactly those subshifts that are a
morphic image of a Sturmian subshift; compare \cite{C,C2,P}.

There are a large number of equivalent characterizations of
Sturmian subshifts; compare \cite{b2}. We are mainly interested in
their geometric description in terms of an irrational rotation.
Let $\alpha \in (0,1)$ be irrational and consider the rotation by
$\alpha$ on the circle,
$$
R_\alpha : [0,1) \to [0,1), \;\; R_\alpha \theta = \{\theta +
\alpha\},
$$
where $\{x\}$ denotes the fractional part of $x$, $\{x\} = x \mod
1$. The coding of the rotation $R_\alpha$ according to a partition
of the circle into two half-open intervals of length $\alpha$ and
$1-\alpha$, respectively, is given by the sequences
$$
v_n(\alpha,\theta) = \chi_{[0,\alpha)}(R_\alpha^n \theta).
$$
We obtain a subshift
\begin{align*}
\Omega_\alpha & = \overline{ \{ v(\alpha,\theta) : \theta \in [0,1) \} }\\
& = \{ v(\alpha,\theta) : \theta \in [0,1) \} \cup \{
\tilde{v}^{(k)}(\alpha) : k \in \Z \} \subset\{0,1\}^{\Z}
\end{align*}
which can be shown to be Sturmian. Here,
$\tilde{v}^{(k)}_n(\alpha) = \chi_{(0,\alpha]}(R_\alpha^{n+k} 0)$.
Conversely, every Sturmian subshift is essentially of this form,
that is, if $\Omega$ is minimal and has complexity function $p(n)
= n+1$, then up to a one-to-one morphism, $\Omega = \Omega_\alpha$
for some irrational $\alpha \in (0,1)$.

By uniform distribution, the frequencies of factors of $\Omega$
are given by the Lebesgue measure of certain subsets of the torus.
Explicitly, if we write $I_0 = [0,\alpha)$ and $I_1 = [\alpha,1)$,
then the word $w = w_1 \ldots w_n \in \{0,1\}^n$ occurs in
$v(\alpha,\theta)$ at site $k+1$ if and only if
$$
\{k\alpha + \theta\} \in I(w_1,\ldots,w_n) := \bigcap_{j = 1}^{n}
R_\alpha^{-j}(I_{w_{j}}).
$$
This shows that the frequency of $w$ is $\theta$-independent and
equal to the Lebesgue measure of $I(w_1,\ldots,w_n)$. It is not
hard to see that $I(w_1,\ldots,w_n)$ is an interval whose boundary
points are elements of the set
$$
P_n(\alpha) := \{ \{ -j \alpha \} : 0 \le j \le n \}.
$$
The $n+1$ points of $P_n(\alpha)$ partition the torus into $n+1$
subintervals and hence the length $h_n(\alpha)$ of the smallest of
these intervals bounds the frequency of a factor of length $n$
from below. It is therefore of interest to study $\limsup n
h_n(\alpha)$.

To this end we recall the notion of a continued fraction
expansion; compare \cite{Khin,RS}. For every irrational $\alpha
\in (0,1)$, there are uniquely determined $a_k \in \N$ such that
\begin{equation}\label{cfe}
\alpha = [a_1,a_2,a_3,\ldots] := \cfrac{1}{a_1+ \cfrac{1}{a_2+
\cfrac{1}{a_3 + \cdots}}} \; .
\end{equation}
The associated rational approximants $\frac{p_k}{q_k}$ are defined
by
\begin{alignat*}{3}
p_0 &= 0, &\quad        p_1 &= 1,   &\quad      p_k &= a_k p_{k-1} + p_{k-2},\\
q_0 &= 1, &             q_1 &= a_1, &           q_k &= a_k q_{k-1}
+ q_{k-2}.
\end{alignat*}
These rational numbers are best approximants to $\alpha$ in the
following sense,
\begin{equation}\label{cfeapp1}
\min_{p,q \in \N \atop 0 < q < q_{k+1}} \left| q \alpha - p
\right| = \left| q_k \alpha - p_k \right|,
\end{equation}
and the quality of the approximation can be estimated according to
\begin{equation}\label{cfeapp2}
\frac{1}{q_k + q_{k+1}} < \left| q_k \alpha - p_k \right| <
\frac{1}{q_{k+1}}.
\end{equation}

By definition, we have
$$
h_n(\alpha) = \min_{0 < |q| \le n} \{ q \alpha \}.
$$
Notice that for $0 < q \le n$, we have $\min \{ \{ q \alpha \} ,
\{ -q \alpha \} \} = \| q\alpha \|$, where we denote $\|x\| =
\min_{p \in \Z} |x - p|$.

In particular,
\begin{equation}\label{htildeh}
h_n(\alpha) = \min_{0 < q \le n} \| q \alpha \|.
\end{equation}
As noted by Hartman \cite{Har}, this shows that $h_n(\alpha)$ can
be expressed in terms of the continued fraction approximants.
Indeed, if we combine \eqref{cfeapp1} and \eqref{htildeh}, we
obtain:

\begin{lemma}[Hartman]\label{hartmanlem}
If $q_k \le n < q_{k+1}$, then
$$
h_n(\alpha) = |q_k \alpha - p_k|.
$$
\end{lemma}

This allows us to show the following:

\begin{theorem}\label{sturmbosh}
Every Sturmian subshift obeys the Boshernitzan condition {\rm
(B)}.
\end{theorem}

\begin{proof}
We only need to show that
\begin{equation}\label{sturmb}
\limsup_{n \to \infty} n h_n(\alpha) \ge C > 0.
\end{equation}
We shall verify this on the subsequence $n_k = q_{k+1} - 1$.
Hartman's lemma together with \eqref{cfeapp2} shows that
$$
n_k h_{n_k}(\alpha) = (q_{k+1} - 1) |q_k \alpha - p_k| \ge
\frac{q_{k+1} - 1}{q_k + q_{k+1}} = \frac{1 - q_{k+1}^{-1}}{1 +
q_k q_{k+1}^{-1}}.
$$
Thus \eqref{sturmb} holds (with $C = 1/3$, say).
\end{proof}

\begin{coro}\label{quasisturmbosh}
Every quasi-Sturmian subshift obeys {\rm (B)}.
\end{coro}

\begin{proof}
This follows from Theorem~\ref{sturmbosh} along with the stability
result, Theorem~\ref{StabilityB}.
\end{proof}

\subsection{Interval Exchange Transformations}\label{IET}

Subshifts generated by interval exchange transformations (IET's)
are natural generalizations of Sturmian subshifts. They were
studied, for example, in
\cite{Bosh3,FHZ,FHZ2,FHZ3,Ke1,Ke2,KN,Ma,R,Vee1,Vee4,Vee2}.

IET's are defined as follows. Given a probability vector $\lambda
= (\lambda_1,\ldots,\lambda_m)$ with $\lambda_i > 0$ for $1 \le i
\le m$, we let $\mu_0 = 0$, $\mu_i = \sum_{j = 1}^i \lambda_j$,
and $I_i = [\mu_{i-1},\mu_i)$. Let $\tau$ be a permutation of
$\CalA_m = \{1,\ldots,m\}$, that is, $\tau \in S_m$, the symmetric
group. Then $\lambda^\tau = (\lambda_{\tau^{-1}(1)}, \ldots,
\lambda_{\tau^{-1}(m)})$ is also a probability vector and we can
form the corresponding $\mu_i^\tau$ and $I_i^\tau$. Denote the
unit interval $[0,1)$ by $I$. The $(\lambda,\tau)$ interval
exchange transformation is then defined by
$$
T : I \to I, \; \; T(x) = x - \mu_{i-1} + \mu_{\tau(i) - 1}^\tau
\text{ for } x \in I_i, \; 1 \le i \le m.
$$
It exchanges the intervals $I_i$ according to the permutation
$\tau$.

The transformation $T$ is invertible and its inverse is given by
the $(\lambda^\tau,\tau^{-1})$ interval exchange transformation.

The symbolic coding of $x \in I$ is $\omega_n(x) = i$ if $T^n(x)
\in I_i$. This induces a subshift over the alphabet $\CalA_m$:
$\Omega_{\lambda,\tau} = \overline{ \{ \omega(x) : x \in I \} }$.

Sturmian subshifts correspond to the case of two intervals as a
first return map construction shows.

Keane \cite{Ke1} proved that if the orbits of the discontinuities
$\mu_i$ of $T$ are all infinite and pairwise distinct, then $T$ is
minimal. In this case, the coding is one-to-one and the subshift
is minimal and aperiodic.  This holds in particular if $\tau$ is
irreducible and $\lambda$ is irrational. Here, $\tau$ is called
irreducible if $\tau(\{1,\ldots,k\}) \not= (\{1,\ldots,k\})$ for
every $k < m$ and $\lambda$ is called irrational if the
$\lambda_i$ are rationally independent.

Regarding property (B), Boshernitzan has proved two results.
First, in \cite{Bosh4} the following is shown:

\begin{theorem}[Boshernitzan]\label{Boshiet1}
For every irreducible $\tau \in S_m$ and for Lebesgue almost every
$\lambda$, the subshift $\Omega_{\lambda,\tau}$ satisfies {\rm
(B)}.
\end{theorem}

In fact, Boshernitzan shows that for every irreducible $\tau \in
S_m$ and for Lebesgue almost every $\lambda$, the subshift
$\Omega_{\lambda,\tau}$ satisfies a stronger condition where the
sequence of $n$'s for which $\eta(n)$ is large cannot be
too sparse. This condition is easily seen to imply (B), and hence
the theorem above.

Note that when combined with Keane's minimality result,
Theorem~\ref{Boshiet1} implies that almost every subshift arising
from an interval exchange transformation is uniquely ergodic. The
latter statement confirms a conjecture of Keane \cite{Ke1} and had
earlier been proven by different methods by Masur \cite{Ma} and
Veech \cite{Vee4}. Keane had in fact conjectured that all minimal
interval exchange transformations would give rise to a uniquely
ergodic system. This was disproved by Keynes and Newton \cite{KN}
using five intervals, and then by Keane \cite{Ke2} using four
intervals (the smallest possible number). The conjecture was
therefore modified in \cite{Ke2} and then ultimately proven by
Masur and Veech.

In a different paper, \cite{Bosh3}, Boshernitzan singles out an
explicit class of subshifts arising from interval exchange
transformations that satisfy (B). The transformation $T$ is said
to be of (rational) rank $k$ if the $\mu_i$ span a $k$-dimensional
space over $\Q$ (the field of rational numbers).

\begin{theorem}[Boshernitzan]
If $T$ has rank $2$, the subshift $\Omega_{\lambda,\tau}$
satisfies {\rm (B)}.
\end{theorem}

\section{Circle Maps}\label{Circle}

Let $\alpha \in (0,1)$ be irrational and $\beta \in (0,1)$
arbitrary. The coding of the rotation $R_\alpha$ according to a
partition into two half-open intervals of length $\beta$ and
$1-\beta$, respectively, is given by the sequences
$$
v_n(\alpha,\beta,\theta) = \chi_{[0,\beta)}(R_\alpha^n \theta).
$$
We obtain a subshift
\begin{equation}\label{cmsubshift}
\Omega_{\alpha,\beta} = \overline{ \{ v(\alpha,\beta,\theta) :
\theta \in [0,1) \} } \subset\{0,1\}^{\Z}.
\end{equation}
Subshifts generated this way are usually called circle map
subshifts or subshifts generated by the coding of a rotation.
These natural generalizations of Sturmian subshifts were studied,
for example, in \cite{AD,AB,BV,DP,Di,Di2,HKS,Kam,Ro}.

To the best of our knowledge, the Boshernitzan condition for this
class of subshifts has not been studied explicitly. It is,
however, intimately related to classical results on inhomogeneous
diophantine approximation problems. In this section we make this
connection explicit and study the condition (B) for circle map
subshifts.

To describe the relation of frequencies of finite words occurring
in a subshift to the length of intervals on the circle, let us
write, in analogy to the Sturmian case, $I_0 = [0,\beta)$ and $I_1
= [\beta,1)$. The word $w = w_1 \ldots w_n \in \{0,1\}^n$ occurs
in $v(\alpha,\beta,\theta)$ at site $k+1$ if and only if
$$
R_\alpha^k (\theta) \in I(w_1,\ldots,w_n) := \bigcap_{j = 1}^{n}
R_\alpha^{-j}(I_{w_{j}}).
$$
Thus the frequency of $w$ is $\theta$-independent and equal to the
Lebesgue measure of $I(w_1,\ldots,w_n)$. Moreover, $I(w_1,\ldots,w_n)$ is an interval whose boundary points are
elements of the set
$$
P_n(\alpha,\beta) := \{ \{ - j\alpha + k\beta \} : 1 \le j \le n,
\; 0 \le k \le 1 \}.
$$
This shows in particular that $\Omega_{\alpha,\beta}$ is
quasi-Sturmian when $\beta \in \Z + \alpha \Z$ as in this case  $P_n(\alpha,\beta)$ splits the  unit interval into  $n+k$ subintervals for large $n$. On the other hand,
when $\beta \not\in \Z + \alpha \Z$, $P_n(\alpha,\beta)$ contains
$2n$ elements and the complexity of $\Omega_{\alpha,\beta}$ is
$p(n) = 2n$ for $n$ large enough.

Again, the points of $P_n(\alpha,\beta)$ partition the torus into
$2n$ (resp., $n + k$) subintervals and hence the length
$h_n(\alpha,\beta)$ of the smallest of these intervals bounds the
frequency of a factor of length $n$ from below. Explicitly, we
have
$$
h_n (\alpha,\beta) = \min \left\{ \left\| q \alpha + r \beta
\right\| : 0 \le |q| \le n, \; 0 \le r \le 1, \; (q,r) \not= (0,0)
\right\}.
$$

Let us also define
$$
\tilde{h}_n(\alpha,\beta) = \min \left\{ \left\| q \alpha + \beta
\right\| : 0 \le |q| \le n \right\}.
$$
Then $h_n(\alpha,\beta) \le \tilde{h}_n(\alpha,\beta)$ and
therefore
\begin{equation}\label{htildehlem}
\limsup_{n \to \infty} n \tilde{h}_n (\alpha,\beta) = 0
\Rightarrow \limsup_{n \to \infty} n h_n (\alpha,\beta) = 0.
\end{equation}

Since we saw in Theorem~\ref{sturmbosh} above that the points of
$P_n(\alpha)$ are nicely spaced for many values of $n$, the
Boshernitzan condition can only fail for a circle map subshift
$\Omega_{\alpha,\beta}$ if the orbit of the $\alpha$-rotation
comes too close to $\beta$. In other words, to prove such a
negative result for a circle map subshift, it should be sufficient
to study $\tilde{h}_n (\alpha,\beta)$, followed by an application
of \eqref{htildehlem}.

Motivated by Hardy and Littlewood \cite{HL}, Morimoto
\cite{Mo,Mo2} carried out an in-depth analysis of the asymptotic
behavior of the numbers $\tilde{h}_n(\alpha,\beta)$. Morimoto's
results and related ones were summarized in \cite{Kok}. While it
is possible to deduce consequences regarding the Boshernitzan
condition from these papers, we choose to give direct and
elementary proofs of our positive results below and make reference
to a specific theorem of Morimoto only for a complementary
negative result.

Our first result shows that the Boshernitzan condition holds in
almost all cases.

\begin{theorem}\label{cmpos}
Let $\alpha \in (0,1)$ be irrational. Then the subshift
$\Omega_{\alpha,\beta}$ satisfies {\rm (B)} for Lebesgue almost
every $\beta \in (0,1)$.
\end{theorem}

\begin{proof}
Denote the set of $\beta$'s for which the Boshernitzan condition
fails by $N(\alpha)$,
$$
N(\alpha) = \{ \beta \in (0,1) : \Omega_{\alpha,\beta} \text{ does
not satisfy (B)} \}.
$$
By \eqref{sturmb} and Theorem~\ref{sturmbosh}, there exists a
sequence $n_k \to \infty$ such that
$$
\liminf_{k \to \infty} n_k {h}_{n_k} (\alpha) = C > 0.
$$
Let $\varepsilon > 0$ with $\varepsilon < C$ be given  and denote the
$\frac{\varepsilon}{2n}$-neighborhood of the set $\{ \{q\alpha\} :
0 < |q| \le n \}$ by $U(\varepsilon,n)$. Clearly, every $\beta \in
N(\alpha)$ belongs to $U(\varepsilon,n_k)$ for $k \ge k_0(\beta)$.
Therefore,
\begin{equation}\label{nalphasub}
N(\alpha) \subseteq \liminf_{k \to \infty} U(\varepsilon,n_k) =
\bigcup_{m=1}^\infty \bigcap_{k \ge m} U(\varepsilon,n_k).
\end{equation}
The sets
$$
S_m = \bigcap_{k \ge m} U(\varepsilon,n_k)
$$
obey $S_m \subseteq S_{m+1}$ and $|S_m| \le \varepsilon$ for every
$m$; $| \cdot |$ denoting Lebesgue measure. Hence,
$$
\left| \liminf_{k \to \infty} U(\varepsilon,n_k) \right| \le
\varepsilon.
$$
It follows that $N(\alpha)$ has zero Lebesgue measure.
\end{proof}

The next result concerns a subclass of $\alpha$'s for which the
Boshernitzan condition holds for all $\beta$'s.

\begin{theorem}
Let $\alpha \in (0,1)$ be irrational with bounded continued
fraction coefficients, that is, $a_n \le C$. Then,
$\Omega_{\alpha,\beta}$ satisfies {\rm (B)} for every $\beta \in
(0,1)$.
\end{theorem}

\begin{proof}
By Lemma~\ref{hartmanlem} and \eqref{cfeapp2}, we have
$$
 {h}_n(\alpha) > \frac{1}{q_k + q_{k+1}},
$$
where $k$ is chosen such that $q_k \le n < q_{k+1}$. Thus, for
every $n$, we have
\begin{equation}\label{bcfproof}
n {h}_n(\alpha) > \frac{n}{q_k + q_{k+1}} \ge
\frac{q_k}{(a_{k+1} + 2) q_k} \ge \frac{1}{C+2}.
\end{equation}

Now assume there exists $\beta \in (0,1)$ such that
$\Omega_{\alpha,\beta}$ does not satisfy {\rm (B)}.  Let
$\varepsilon = (7C+14)^{-1}$. As $\limsup_{n\to \infty} n
\varepsilon (n) =0$, we have $n \varepsilon (n) < \epsilon$ for
every sufficiently large $n$. Thus, for each such $n$ we can find
a word of length $n$ with frequency less than $\varepsilon/n $.
Now, each such word corresponds to an interval with length less
than $\varepsilon /n $ with boundary points in $P_n (\alpha,
\beta)$. Moreover, invoking \eqref{bcfproof} and the fact that
$\epsilon < 1 / (C +  2 )$, we infer that the length of the
interval has the form $|m_n \alpha - \beta - k_n| $ with
$|m_n|\leq n$.  To summarize, we see that for every $n$ large
enough there exist $k_n, m_n$ with $|m_n| \le n$ such that
$$
|m_n \alpha - \beta - k_n| \le \frac{\varepsilon}{n}.
$$
Clearly, the mapping $n \mapsto
m_n$ can take on each value only finitely many times. Therefore,
there exists a sequence $n_j \to \infty$ such that $m_{n_j} \not=
m_{n_j + 1}$. This implies
\begin{align*}
\left| \left( m_{n_j + 1} - m_{n_j} \right) \alpha - \left( k_{n_j
+ 1} - k_{n_j} \right) \right| & \le \left| m_{n_j + 1} \alpha -
\beta - k_{n_j + 1} \right| + \left| m_{n_j} \alpha - \beta -
k_{n_j} \right|\\
& \le \frac{\varepsilon}{n_j + 1} + \frac{\varepsilon}{n_j}\\
& \le \frac{2\varepsilon}{n_j}.
\end{align*}
Since $0 < \left| m_{n_j + 1} - m_{n_j} \right| \le 2(n_j + 1) \le
3 n_j=:\tilde{n}_j$, we obtain $\tilde{n}_j h_{\tilde{n}_j}
(\alpha) \le 6 \varepsilon < (C+2)^{-1}$, which contradicts
\eqref{bcfproof}.
\end{proof}

This raises the question whether $\Omega_{\alpha,\beta}$ satisfies
(B) for every $\beta$ also in the case where $\alpha$ has
unbounded coefficients $a_n$. It is a consequence of a result of
Morimoto \cite{Mo2} that this is not the case.

\begin{theorem}[Morimoto]\label{morimoto}
Let $\alpha \in (0,1)$ be irrational with unbounded continued
fraction coefficients. Then, there exists $\beta \in (0,1)$ such
that
$$
\limsup_{n \to \infty} n \tilde{h}_n (\alpha,\beta) = 0.
$$
\end{theorem}

\begin{coro}\label{cmneg}
Let $\alpha \in (0,1)$ be irrational with unbounded continued
fraction coefficients. Then, there exists $\beta \in (0,1)$ such
that $\Omega_{\alpha,\beta}$ does not satisfy {\rm (B)}.
\end{coro}

\begin{proof}
This is an immediate consequence of Theorem~\ref{morimoto} and
\eqref{htildehlem}.
\end{proof}

We close this section with a brief discussion of the case where
the circle is partitioned into a finite number of half-open
intervals. To be specific, let $0 < \beta_1 < \cdots < \beta_{p-1}
< 1$ and associate the intervals of the induced partition with $p$
symbols: Let $\beta_p = \beta_0 = 0$ and
$$
v_n(\theta) = k \Leftrightarrow R_\alpha^n (\theta) \in
[\beta_k,\beta_{k+1}).
$$
We obtain a subshift over the alphabet $\{0,1,\ldots, p-1\}$,
$$
\Omega_\beta = \overline{ \{ v(\theta) : \theta \in [0,1) \} }.
$$
Again, the word $w = w_1 \ldots w_n \in \{0,1\}^n$ occurs in
$v(\theta)$ at site $k+1$ if and only if
$$
R_\alpha^k \theta \in I(w_1,\ldots,w_n) := \bigcap_{j = 1}^{n}
R_\alpha^{-j}(I_{w_{j}})
$$
and the connected components of the sets $I(w_1,\ldots,w_n)$ are
bounded by the points
\begin{equation}\label{multiintpoint}
\{ - j \alpha + \beta_k : 1 \le j \le n, 0 \le k \le p-1 \}.
\end{equation}

Recall that $\limsup_{n \to \infty} n \tilde{h}_n (\alpha,\beta)$
is an important quantity in the case of a partition of the circle
into two intervals. In fact, we showed that this quantity being
positive is a necessary condition for (B) to hold. When there are
three or more intervals, however, we will need to require a much
stronger condition as the $\beta_i$'s may now ``take turns'' in
being well approximated by the $\alpha$-orbit. Indeed, we shall
now be interested in studying $\liminf_{n \to \infty} n
\tilde{h}_n (\alpha,\gamma)$ (for certain values of $\gamma$,
associated with the $\beta_i$'s). More precisely, define the
following quantity:
$$
M(\alpha,\gamma) = \liminf_{|n| \to \infty} |n| \cdot \| n\alpha -
\gamma\|.
$$
Let
$$
P(\alpha) = \{ \gamma : M(\alpha,\gamma) > 0\}.
$$
Then, we have the following result:

\begin{theorem}
Let $\alpha \in (0,1)$ be irrational. Suppose that $0 = \beta_0 <
\beta_1 < \cdots < \beta_{p-1} < \beta_p = 1$ are such that
$$
\beta_k - \beta_l \in P(\alpha) \text{ for } 0 \le k \not= l \le
p-1.
$$
Then the subshift $(\Omega_\beta,T)$ satisfies the Boshernitzan
condition {\rm (B)}.
\end{theorem}

\noindent\textit{Remarks.} 1. This gives a finite number of
conditions whose combination is a sufficient condition for (B) to
hold.
\\[1mm]
2. The set $P(\alpha)$ is non-empty for every irrational $\alpha$.
In fact, for every irrational $\alpha$ there exists a suitable
$\gamma$ such that $M(\alpha,\gamma) > 1/32$; compare
\cite[Theorem~IV.9.3]{RS}.
\\[1mm]
3. We discuss in Appendix~\ref{pinnerapp} how $M(\alpha,\gamma)$
can be computed with the help of the so-called \textit{negative
continued fraction expansion} of $\alpha$ and the
$\alpha$-\textit{expansion} of $\gamma$.

\begin{proof}
By \eqref{multiintpoint}, all frequencies of words of length $n$
are bounded from below by
$$
\hat{h}_n(\alpha,\beta) = \min \left\{ \left\| q \alpha + \beta_k
- \beta_l \right\| : 0 \le |q| \le n, \; 0 \le k,l \le p-1, \;
(q,k-l) \not= (0,0) \right\}.
$$
As in our considerations above, we choose a sequence $n_k \to
\infty$ such that
$$
\liminf_{k \to \infty} n_k {h}_{n_k} (\alpha) = C > 0.
$$
By assumption, we have
$$
D = \min \{ M(\alpha,\beta_k - \beta_l) : 0 \le k \not= l \le p-1
\}
> 0.
$$
Notice that with these choices of $C$ and $D$, frequencies of
words of length $n_k$ are bounded from below by
$$
\hat{h}_{n_k} (\alpha,\beta) \ge \min \left\{ \frac{C -
o(1)}{n_k}, \frac{D - o(1)}{n_k} \right\}.
$$
Putting everything together, we obtain
$$
\limsup_{n \to \infty} n \cdot \eta(n) \ge \liminf_{k \to
\infty} n_k \cdot \eta(n_k) \ge \min \{C,D\} > 0,
$$
and hence (B) is satisfied.
\end{proof}

\section{Arnoux-Rauzy Subshifts and Episturmian Subshifts}\label{Arnoux-Rauzy}

In this section we consider another natural generalization of
Sturmian subshifts, namely, Arnoux-Rauzy subshifts and, more
generally, episturmian subshifts. These subshifts were studied,
for example, in \cite{AR,DZ,DJP,JP,JV,RZ,WZ}. They share with
Sturmian subshifts the fact that, for each $n$, there is a unique
subword of length $n$ that has multiple extensions to the right.
Our main results will show that, similarly to the circle map case,
the Boshernitzan condition is almost always satisfied, but not always.

Let us consider a minimal subshift $\Oomega$ over the alphabet
$\CalA_m = \{ 1,2,\ldots,m\}$, where $m \ge 2$. Recall that the
set of subwords of length $n$ occurring in elements of $\Omega$ is
denoted by $\CalW_n (\Omega)$ (cf.~\eqref{subwordsomega}) and that
the complexity function $p$ is defined by $p(n) = \sharp \CalW_n
(\Omega) $ (cf.~\eqref{complexitydef}). A word $w \in \CalW
(\Omega)$ is called \textit{right-special} (resp.,
\textit{left-special}) if there are distinct symbols $a,b \in
\CalA_m$ such that $wa,wb \in \CalW (\Omega)$ (resp., $aw,bw \in
\CalW (\Omega)$). A word that is both right-special and
left-special is called \textit{bispecial}.

For later use, let us recall the \textit{Rauzy graphs} that are
associated with $\CalW (\Omega)$. For each $n$, we consider the
directed graph $\mathcal{R}_n = (V_n, A_n)$, where the vertex set
is given by $V_n = \CalW_n (\Omega)$, and $A_n$ contains the arc
from $aw$ to $wb$, $a,b \in \CalA_m$, $|w| = n-1$, if and only if
$awb \in \CalW_{n+1} (\Omega)$. That is, $|V_n| = p(n)$ and $|A_n|
= p(n+1)$. Moreover, a word is right-special (resp., left-special)
if and only if its out-degree (resp., in-degree) is $\ge 2$.

Note that the complexity function of a Sturmian subshift obeys
$p(n+1) - p(n) = 1$ for every $n$ and hence for every length,
there is a unique right-special factor and a unique left-special
factor, each having exactly two extensions. This property is
clearly characteristic for a Sturmian subshift.

Arnoux-Rauzy subshifts and episturmian subshifts relax this
restriction on the possible extensions somewhat, and they are
defined as follows: $\Omega$ is called an \textit{Arnoux-Rauzy
subshift} if for every $n$, there is a unique right-special word
$r_n \in \CalW (\Omega)$ and a unique left-special word $l_n \in
\CalW (\Omega)$, both having exactly $m$ extensions. This implies
in particular that $p(1) = m$ and hence
$$
p(n) = (m-1)n + 1.
$$
Arnoux-Rauzy subshifts over $\CalA_2$ are exactly the Sturmian
subshifts.

On the other hand, $\Omega$ is called \textit{episturmian} if
$\CalW (\Omega)$ is closed under reversal (i.e., for every $w =
w_1 \ldots w_n \in \CalW (\Omega)$, we have $w^R = w_n \ldots w_1
\in \CalW (\Omega)$) and for every $n$, there is exactly one
right-special word $r_n \in \CalW (\Omega)$.

It is easy to see that every Arnoux-Rauzy subshift is episturmian.
On the other hand, every episturmian subshift turns out to be a
morphic image of some Arnoux-Rauzy subshift. We shall explain this
connection below. Since we are interested in studying the
Boshernitzan condition, this fact is important and allows us to
limit our attention to the Arnoux-Rauzy case.

Risley and Zamboni \cite{RZ} found two useful descriptions of a
given Arnoux-Rauzy subshift, namely, in terms of the recursive
structure of the bispecial words and in terms of an $S$-adic
system.

Let $\epsilon$ be the empty word and let  $\{ \epsilon = w_1, w_2,
\ldots \}$ be the set of all bispecial words in $\CalW (\Omega)$,
ordered so that $0 = |w_1| < |w_2| < \cdots$. Let $I = \{i_n\}$ be
the sequence of elements $i_n$ of $\CalA_m$ so that $w_n i_n$ is
left-special. The sequence $I$ is called the \textit{index
sequence} associated with $\Omega$. Risley and Zamboni prove that,
for every $n$, $w_{n+1}$ is the \textit{palindromic closure} $(w_n
i_n)^+$ of $w_n i_n$, that is, the shortest palindrome that has
$w_n i_n$ as a prefix. Conversely, given any sequence $I$, one can
associate a subshift $\Omega$ as follows: Start with $w_1 =
\epsilon$ and define $w_n$ inductively by $w_{n+1} = (w_n i_n)^+$.
The sequence of words $\{ w_n \}$ has a unique one-sided infinite
limit $w_{\infty} \in \CalA_m^{\N}$, called the
\textit{characteristic sequence}, which then gives rise to the
subshift $(\Omega(I),T)$ in the standard way; $\Omega(I)$ consists
of all two-sided infinite sequences whose subwords occur in $w$.
Risley and Zamboni prove the following characterization.

\begin{prop}[Risley-Zamboni]\label{arbisp}
For every Arnoux-Rauzy subshift $\Oomega$ over $\CalA_m$, every $a
\in \CalA_m$ occurs in the index sequence $\{i_n\}$ infinitely
many times and $\Omega = \Omega(I)$. Conversely, for every
sequence $\{i_n\} \in \CalA_m^{\N}$ such that every $a \in
\CalA_m$ occurs in $\{i_n\}$ infinitely many times,
$(\Omega(I),T)$ is an Arnoux-Rauzy subshift and $\{i_n\}$ is its
index sequence.
\end{prop}

The $S$-adic description of an Arnoux-Rauzy subshift, that is,
involving iterated morphisms chosen from a finite set, found in
\cite{RZ} reads as follows.

\begin{prop}[Risley-Zamboni]\label{aradic}
Let $\Oomega$ be an Arnoux-Rauzy subshift over $\CalA_m$ and
$\{i_n\}$ the associated index sequence. For each $a \in \CalA_m$,
define the morphism $\tau_a$ by
$$
\tau_a(a) = a \text{ and } \tau_a(b) = ab \text{ for } b \in
\CalA_m \setminus \{a\}.
$$
Then for every $a \in \CalA_m$, the characteristic sequence is
given by
$$
\lim_{m \to \infty} \tau_{i_1} \circ \cdots \circ \tau_{i_m} (a).
$$
\end{prop}

We can now state our positive result regarding the Boshernitzan
condition for Arnoux-Rauzy subshifts.

\begin{theorem}\label{arpos}
Let $\Oomega$ be an Arnoux-Rauzy subshift over $\CalA_m$ and
$\{i_n\}$ the associated index sequence. Suppose there is $N \in
\N$ such that for a sequence $k_j \to \infty$, each of the words
$i_{k_j} \ldots i_{k_j + N - 1}$ contains all symbols from
$\CalA_m$. Then the Boshernitzan condition {\rm (B)} holds.
\end{theorem}

This result is similar to Theorem~\ref{cmpos} in the sense that if
we put any probability measure $\nu$ on $\CalA_m$ assigning
positive weight to each symbol, then almost all sequences
$\{i_n\}$ with respect to the product measure $\nu^{\N}$
correspond to Arnoux-Rauzy subshifts that satisfy the assumption
of Theorem~\ref{arpos}.

Before proving this theorem, we state our negative result, which
is an analog of Corollary~\ref{cmneg}.

\begin{theorem}\label{arneg}
For every $m \ge 3$, there exists an Arnoux-Rauzy subshift over
$\CalA_m$ that does not satisfy the Boshernitzan condition {\rm
(B)}.
\end{theorem}

\begin{remark}{\rm The assumption $m \ge 3$ is of course necessary since the case $m
= 1$ is trivial and the case $m =2$ corresponds to the Sturmian
case, where the Boshernitzan condition always holds; compare
Theorem~\ref{sturmbosh}.}
\end{remark}

The Arnoux-Rauzy subshifts are uniquely ergodic and we set
$$ d(w) \equiv \mu (V_w), \, w\in \CalW (\Omega),$$
where, as usual, the unique invariant probability measure is denoted by $\mu$.

\begin{proof}[Proof of Theorem~\ref{arpos}.]
This proof employs the description of the subshift in terms of the
bispecial words; compare Proposition~\ref{arbisp}.

Observe that there is some $k_0$ such that $|w_k| \le 2 |w_{k-1}|$
for every $k \ge k_0$. Essentially, we need that $i_1, \ldots,
i_{k_0 - 1}$ contains all symbols from $\CalA_m$.

Now consider a value of $k \ge k_0$ such that $i_{k} \ldots
i_{k+N-1}$ contains all symbols from $\CalA_m$. We claim that
\begin{equation}\label{arposgoal}
|w_k| \cdot \eta (|w_k|) \ge 2^{-N}.
\end{equation}
By the assumption, this implies
$$
\limsup_{n \to \infty} n \cdot \eta (n) \ge 2^{-N}
$$
and hence the Boshernitzan condition (B).

The Rauzy graph $\mathcal{R}_{|w_k|}$ has one vertex (namely,
$w_k$) with in-degree and out-degree $m$, while all other vertices
have in-degree and out-degree $1$. Thus, the graph splits up into
$m$ loops that all contain $w_k$ and are pairwise disjoint
otherwise. These loops can be indexed in an obvious way by the
elements of the alphabet $\CalA_m$.

Since $w_{k+1} = (w_k i_k)^+$, $w_{k+1}$ begins and ends with
$w_k$ and, moreover, $w_{k+1}$ contains all words that correspond
to the loop in $\mathcal{R}_{|w_k|}$ indexed by $i_k$. Iterating
this argument, we see that $w_{k + N}$ contains the words from all
loops and hence all words from $\CalW_{|w_k|} (\Omega)$. This
implies
$$
\min_{w \in \CalW_{|w_k|} (\Omega)} d(w) \ge d(w_{k+N}) \ge
\frac{1}{|w_{k+N}|} \ge \frac{1}{2^N |w_k|}
$$
and hence \eqref{arposgoal}, finishing the proof.
\end{proof}

\begin{proof}[Proof of Theorem~\ref{arneg}.]
This proof employs the description of the subshift in terms of an
$S$-adic structure; compare Proposition~\ref{aradic}.

We shall construct an index sequence $\{i_n\}$ over three symbols
(i.e., over the alphabet $\CalA_3$) such that the corresponding
Arnoux-Rauzy subshift does not satisfy the Boshernitzan condition
(B). It is easy to verify that the same idea can be used to
construct such a subshift over $\CalA_m$ for any $m \ge 3$.

The index sequence will have the form
\begin{equation}\label{indexsequ}
i_1 i_2 i_3 \ldots = 1^{a_1} 2^{a_2} 3^{a_3} 1^{a_4} 2^{a_5}
3^{a_6} 1^{a_7} \ldots
\end{equation}
with a rapidly increasing sequence of integers, $\{a_n\}$.

By the special form of the Rauzy graph, the words $w_k a$ label
all the frequencies of words in $\CalW_{|w_k|+1} (\Omega)$ since
words corresponding to arcs on a given loop in
$\mathcal{R}_{|w_k|}$ must have the same frequency. Put
differently,
\begin{equation}\label{arfreqloop}
\eta( |w_k| + 1 ) = \min_{a \in \CalA_3} d_{w_\infty} (w_k a).
\end{equation}
Here, we make the dependence of the frequency on $w_{\infty}$ explicit.

Moreover, it is sufficient to control $\eta(n)$ for these
special values of $n$ since every subword $u$ that is not
bispecial has a unique extension to either the left or the right,
and this extension must have the same frequency. This shows
\begin{equation}\label{goodsubs}
\eta( |w_k| + 1 ) \ge \eta( n ) \text{ for } |w_k| +
1 \le n \le |w_{k+1}|.
\end{equation}

Now write $\mu_{k,m} = \tau_{i_k} \circ \cdots \circ \tau_{i_{k +
m - 1}}$. Proposition~\ref{aradic} says that the characteristic
sequence is given by the limit
$$
w = \lim_{m \to \infty} \mu_{1,m} (a) \text{ for every } a \in
\CalA_3.
$$
We also define
$$
w^{(k)} = \lim_{m \to \infty} \mu_{k,m} (a) =
(\mu_{1,k-1})^{-1}(w).
$$
By \cite{JP}, $w^{(k)}$ is the derived sequence labeling the
return words of $w_k$ in $w_{\infty}$. In particular, $w^{(k)}$ labels the
occurrences of $w_k a$, $a \in \CalA_3$, in $w$. Moreover,
\begin{equation}\label{arfreqest}
d_{w_\infty} (w_k a) = \frac{d_{w^{(k)}} (a)}{\sum_{b \in \CalA_3}
d_{w^{(k)}} (b) |\mu_{1,k-1}(b)|} \le \frac{d_{w^{(k)}}
(a)}{\min_{b \in \CalA_3} |\mu_{1,k-1}(b)|} .
\end{equation}
Combining \eqref{arfreqloop} and \eqref{arfreqest}, we obtain
\begin{equation}\label{arfreqest2}
( |w_k| + 1 ) \cdot \eta( |w_k| + 1 ) \le \frac{|w_k| +
1}{\min_{b \in \CalA_3} |\mu_{1,k-1}(b)|} \cdot \min_{a \in
\CalA_3} d_{w^{(k)}} (a).
\end{equation}
Notice that $(|w_k| + 1)(\min_{b \in \CalA_3}
|\mu_{1,k-1}(b)|)^{-1}$ only depends on $i_1, \ldots, i_{k-1}$ and
$\min_{a \in \CalA_3} d_{w^{(k)}} (a)$ only depends on $i_k,
i_{k+1}, \ldots$. Thus, if we choose a rapidly increasing sequence
$\{a_n\}$ in \eqref{indexsequ}, we can arrange for
\begin{equation}\label{zeroonsubs}
\lim_{k \to \infty} ( |w_k| + 1 ) \cdot \eta( |w_k| + 1 ) =
0.
\end{equation}
This together with \eqref{goodsubs} implies
$$
\lim_{n \to \infty} n \cdot \eta( n ) = 0,
$$
proving the theorem.

Let us briefly comment on \eqref{zeroonsubs}. Choose a
monotonically decreasing sequence $e_k \to 0$. Assign any value
$\ge 1$ to $a_1$. Then, $a_2$ should be chosen large enough so
that for $1 \le k \le a_1$, \eqref{arfreqest2} yields
\begin{equation}\label{arneeded}
( |w_k| + 1 ) \cdot \eta( |w_k| + 1 ) \le e_k.
\end{equation}
Here we use that between consecutive $3$'s in $w^{(k)}$, there
must be at least $a_2$ $2$'s. Next, we choose $a_3$ so large that
\eqref{arneeded} holds for $a_1 + 1 \le k \le a_2$. Here we use
that between consecutive $1$'s in $w^{(k)}$, there must be at
least $a_3$ $3$'s. We can continue in this fashion, thereby
generating a sequence $\{a_n\}$ such that \eqref{arneeded} holds
for all $k$. This shows in particular that $( |w_k| + 1 ) \cdot
\eta( |w_k| + 1 )$ can go to zero arbitrarily fast.
\end{proof}

One may wonder what sequences are generated by the procedures
described before Propositions~\ref{arbisp} and \ref{aradic} if one
starts with an index sequence that does not necessarily satisfy
the assumption above, namely, that all symbols occur infinitely
often. It was shown by Droubay, Justin and Pirillo \cite{DJP,JP}
that one obtains episturmian subshifts and, conversely, every
episturmian subshift can be generated in this way.

\begin{prop}[Droubay, Justin, Pirillo]\label{epibisp}
For every episturmian subshift $\Oomega$ over $\CalA_m$, there
exists an index sequence $\{i_n\}$ such that $\Omega = \Omega(I)$.
Conversely, for every sequence $\{i_n\} \in \CalA_m^{\N}$,
$(\Omega(I),T)$ is an episturmian subshift and $\{i_n\}$ is its
index sequence. For every $a \in \CalA_m$, the characteristic
sequence is given by
$$
\lim_{m \to \infty} \tau_{i_1} \circ \cdots \circ \tau_{i_m} (a).
$$
\end{prop}

We can now quickly deduce results concerning (B) (and hence (P))
for episturmian subshifts. If $\Oomega$ is an episturmian subshift
over $\CalA_m$, denote by $\CalA \subseteq \CalA_m$ the set of all
symbols that occur in its index sequence infinitely many times.
Fix $k$ such that $i_k, i_{k+1}, \ldots$ only contains symbols
from $\CalA$. Thus this tail sequence corresponds to an
Arnoux-Rauzy subshift over $|\CalA|$ symbols and the given
episturmian subshift is a morphic image (under $\mu_{1,k-1}$) of
it. (Note that $|\CalA| \ge 2$ since $\Oomega$ is aperiodic.) If
the associated Arnoux-Rauzy subshift satisfies (B) (if, e.g.,
Theorem~\ref{arpos} applies), then $\Oomega$ satisfies (B) by
Theorem~\ref{StabilityB}. On the other hand, since every
Arnoux-Rauzy subshift is episturmian, Theorem~\ref{arneg} shows
that not all episturmian subshifts satisfy (B). In this context,
it is interesting to note that Justin and Pirillo showed that all
episturmian subshifts are uniquely ergodic \cite{JP}.

\section{Application to Schr\"odinger Operators} \label{Application}
In this section we discuss applications of our previous study to
spectral theory of Schr\"odinger operators. This is based on
methods introduced in \cite{Len2} by Lenz.

\smallskip

Let $\Oomega$ be a minimal uniquely ergodic subshift over the finite set  $\CalA$ and
assume $\CalA \subset \R$. As discussed in the introduction, $\Oomega$
gives rise to the family $(H_\omega)_{\omega\in \Omega}$ of
selfadjoint operators $H_\omega : \ell^2 (\Z)\longrightarrow \ell^2
(\Z)$ acting by
$$(H_V u )(n) \equiv u(n+1) + u(n-1) + \omega (n) u(n).$$
As $\Oomega$ is minimal, there exists a set $\Sigma\equiv \Sigma(\Oomega) \subset \R$ with
$$\sigma(H_\omega) = \Sigma\;\:\mbox{for all $\omega \in \Omega$}$$
(see, e.g., \cite{BIST}). We will assume furthermore that
$\Oomega$ is aperiodic. Such subshifts have attracted a lot of
attention in recent years for both physical and mathematical
reasons:

These subshifts can serve as models for a special class of solids
discovered in 1984 by Shechtman et al.\ \cite{SBGC}. These solids,
later called quasicrystals, have very special mechanical,
electrical, and diffraction properties \cite{Jan,Sen}. In the
quantum mechanical description of electrical (i.e., conductance)
properties of these solids, one is led to the operators
$(H_\omega)$ above. These operators in turn have a tendency to
display intriguing mathematical features. These features include:

\begin{itemize}

\item[$\Zet$] Cantor spectrum of Lebesgue measure zero, that is,
$\Sigma$ is a Cantor set of Lebesgue measure zero.

\item[$\SingCont$] Purely singular continuous spectrum, that is,
absence of both point spectrum and absolutely continuous spectrum.

\item[$\AnTr$] Anomalous transport.

\end{itemize}

By now, absence of absolutely continuous spectrum is completely
established for all relevant subshifts due to results of Last and
Simon \cite{LS} in combination with earlier results of Kotani
\cite{Kot}. The other spectral features have been investigated for
large, but special, classes of examples.  Here, our focus is on
$\Zet$. As for the other properties, we refer the reader to the
survey articles \cite{Dam,Sut}.

\smallskip

The property $\Zet$ has been investigated for several models: For
the period-doubling substitution and the Thue-Morse substitution,
it was shown to hold by Bellissard et al.\ in \cite{BBG}
(cf.~earlier work of Bellissard \cite{Bel} as well). A more
general result for primitive substitutions has then been obtained
by Bovier and Ghez \cite{BG}. Recently, proofs of $\Zet$ for all
primitive substitutions were obtained by Liu et al.\ \cite{LTWW}
and, independently, by Lenz \cite{Len2}. For special examples of
on-primitive substitutions, $\Zet$ has recently been investigated
by de Oliveira and Lima \cite{OL2}. Their results were extended by
Damanik and Lenz \cite{DL5}.

For Sturmian operators, $\Zet$ has been proven by S\"ut\H{o} in
the golden mean case ($=$ Fibonacci substitution) \cite{Sut,Sut2}.
The general case was then treated by Bellissard et al.\
\cite{BIST}. A different approach to $\Zet$ in the Sturmian case
has been developed in \cite{DL6} by Damanik and Lenz.  A suitably
modified version of this approach can also be used to study $\Zet$
for a certain class of substitutions as shown by Damanik
\cite{Dam6}.

For quasi-Sturmian operators, $\Zet$ was shown in \cite{DL4}. Later a different proof was given in \cite{Len3}.

\smallskip

All approaches to $\Zet$ are based on a fundamental result of Kotani
\cite{Kot}. To discuss this result, we need some preparation.

Spectral properties of the operators $(H_\omega)$ are intimately linked
to behavior of solutions of the difference equation
 \begin{equation}\label{gleichung}
u(n+1) + u(n-1) + (\omega(n) -E) u (n)= 0
\end{equation}
for $E\in \R$.  To study this behavior, we define, for $E\in \R$,
the locally constant function $M^E : \Omega \longrightarrow \SL$ by
\begin{equation}  \label{transfer}
M^E(\omega) \equiv\left( \begin{array}{cc} E-\omega(1)  & -1\\1 & 0 \end{array}   \right).
\end{equation}
Then, it is easy to see that a sequence $u$ is a solution of the
difference equation \eqref{gleichung} if and only if
\begin{equation} \label{wichtig}
\left( \begin{array}{c} u(n+1) \\ u(n) \end{array} \right) = M^E
(n,\omega)\left( \begin{array}{c} u(1) \\ u(0) \end{array} \right), \:
n\in \Z.
\end{equation}
The rate of exponential growth of solutions of \eqref{gleichung}
is then measured by the so-called Lyapunov exponent
$\gamma(E)\equiv \Lambda(M^E)$. The fundamental result of Kotani,
mentioned above, says that (due to aperiodicity)
\begin{equation}\label{Kotani}
|\{E\in \R : \gamma(E) =0\} | = 0,
\end{equation}
where $| \cdot|$ denotes Lebesgue measure on $\R$. By general
principles, it is clear that $\{E\in \R : \gamma(E) =0\} \subset
\Sigma$ \cite{CL}.  The overall strategy to prove $\Zet$ is then to
show
\begin{equation}\label{SigmaundLE}
\Sigma =  \{E\in \R : \gamma(E) = 0\}.
\end{equation}
Given \eqref{SigmaundLE}, $\Sigma$ can not contain
an interval by \eqref{Kotani}. Moreover, $\Sigma$ is a closed set, as
the spectrum of an operator always is. Finally, $\Sigma$ does not
contain isolated points, again by general principles on random
operators \cite{CL}.  Hence, $\Sigma$ is a Cantor set of measure zero
if \eqref{SigmaundLE} holds.

\smallskip

The standard approach to \eqref{SigmaundLE} used to rely on trace
maps.  Trace maps are a powerful tool in the study of spectral
properties. In particular, they can be used not only to study
$\Zet$, but also to investigate $\SingCont$ and $\AnTr$. However,
trace maps do not seem to be available as soon as the dynamical
systems get more complicated.  This difficulty is avoided in a new
approach to $\Zet$ introduced in \cite{Len2}. There, validity of
\eqref{SigmaundLE} is related to certain ergodic properties of the
underlying dynamical system.  More precisely, the abstract
cornerstone of this new approach is the following result.

\begin{theorem}\label{Nullstelle} \cite{Len2}
Let $\Oomega$ be a minimal uniquely ergodic subshift over $\CalA
\subset \R$. Then, $\Sigma = \{E \in \R : \gamma(E) = 0\}$ if and only
if $M^E$ is uniform for every $E\in \R$. In this case, the map $\gamma
: \R \longrightarrow [0,\infty)$ is continuous.
\end{theorem}

\smallskip

Given this theorem, it becomes possible to show \eqref{SigmaundLE}
by studying uniformity of the functions $M^E$. In fact, as shown
in \cite{Len2}, uniformity of $M^E$ holds for all systems
satisfying (PW) and, in particular, for all linearly repetitive
systems (see \cite{Len1} as well).  In \cite{Len2}, this was used
to prove $\Zet$ for all primitive substitutions.  Later $\Zet$ has
been established for various further systems by showing linear
repetitivity \cite{AD,DL5,OL2}.

\begin{proof}[Proof of Theorem \ref{zeromeasure}]
Given Theorem \ref{Nullstelle} and Kotanis result \eqref{Kotani},
the assertion follows easily from our main result: By (B) and
Theorem \ref{main}, the function $M^E$ is uniform for every $E\in
\R$. By Theorem \ref{Nullstelle}, this implies $\Sigma=\{E :
\gamma(E) =0\}$. By \eqref{Kotani}, this gives that $\Sigma$ is a
Cantor set of Lebesgue measure zero, as discussed above.
\end{proof}

\begin{proof}[Proof of Theorem \ref{CantorCircleMap}]
The result follows from Theorem~\ref{zeromeasure} and the results
regarding the validity of (B) for circle map subshifts of
Section~\ref{Circle}.
\end{proof}

\begin{proof}[Proof of Theorem \ref{ContinuityLE}]
The assertion is immediate from Theorem~\ref{Nullstelle} and our
main result, Theorem \ref{main}.
\end{proof}

As becomes clear from the discussion in Section \ref{Known},
Theorem \ref{zeromeasure} generalizes all earlier results on
$\Zet$. Moreover, it gives various new ones. One of these new
results on $\Zet$ is Theorem \ref{CantorCircleMap}. Similarly,
combining Theorem \ref{zeromeasure} and the results of
Subsection~\ref{IET}, we obtain another new result on $\Zet$ for
subshifts associated with interval exchange transformations. To
the best of our knowledge this is the first result on $\Zet$ for
operators associated to interval exchange transformations (not
counting those which are Sturmian or linearly repetitive).

\begin{theorem}
Let $\tau \in S_m$ be irreducible. Then, for Lebesgue almost every
$\lambda$, $\Sigma = \Sigma( \Omega_{\lambda,\tau})$ is a Cantor
set of Lebesgue measure zero.
\end{theorem}

\begin{proof}
As discussed in Subsection~\ref{IET}, if $\tau \in S_m$ is
irreducible, $\Omega_{\lambda,\tau}$ is minimal, aperiodic, and
satisfies (B) for almost every $\lambda$. This, combined with
Theorem~\ref{zeromeasure}, yields the assertion.
\end{proof}

Finally, we also mention the following result for Arnoux-Rauzy
subshifts, which follows from Theorem \ref{zeromeasure} and the
discussion in Section \ref{Arnoux-Rauzy}

\begin{theorem} Let $\Oomega$ be an aperiodic  Arnoux-Rauzy subshift over $\CalA_m$ and
$\{i_n\}$ the associated index sequence. Suppose there is $N \in
\N$ such that for a sequence $k_j \to \infty$, each of the words
$i_{k_j} \ldots i_{k_j + N - 1}$ contains all symbols from
$\CalA_m$. Then, $\Sigma$ is a Cantor set of measure zero.
\end{theorem}

\bigskip

\noindent\textit{Acknowledgments.}  We thank Barry Simon for
stimulating discussions. A substantial part of this work was done
while one of the authors (D.L.) was visiting Caltech in September
2003. He would like to thank Barry Simon and the Department of
Mathematics for the warm hospitality.

\begin{appendix}

\section{Almost Every Circle Map Subshift Has Infinite Index}\label{highpowers}

In this section, we show that the previous results on zero-measure
Cantor spectrum for Schr\"odinger operators associated with circle
map subshifts only cover a zero-measure set in parameter space.
This should be seen in connection with
Theorem~\ref{CantorCircleMap}, where this spectral result is
established for almost all parameter values.

Recall that for every $(\alpha,\beta) \in (0,1) \times (0,1)$, we
may define a subshift $\Omega_{\alpha,\beta}$ as in
\eqref{cmsubshift}.

Proofs of zero-measure spectrum for the associated operators based
on trace map dynamics were given in \cite{BIST,DL4,Sut,Sut2}. They
cover the case of arbitrary irrational $\alpha \in (0,1)$ and
$\beta$'s in $(0,1)$ of the form $\beta = m \alpha + n$. This is
clearly a zero-measure set in $(0,1) \times (0,1)$.

The paper \cite{AD} applies the results of \cite{Len2} and shows
zero-measure spectrum for a class of circle map subshifts that is
characterized by means of a generalized continued fraction
algorithm. Essentially, \cite{AD} characterizes the pairs
$(\alpha,\beta)$ for which the associated subshifts are linearly
recurrent. We want to show that these, too, form a set of measure
zero.

To this end, we note that every aperiodic linearly recurrent
subshift $\Omega$ has finite index in the sense that there is $N <
\infty$ such that its set of finite subwords, $\CalW (\Omega)$,
contains no word of the form $w^N$. (This is immediate from the
definition.)

We say that a subshift $\Omega$ has infinite index if for every $n
\ge 1$, there is a word $w$ such that $w^n \in \CalW (\Omega)$ and
prove the following:

\begin{prop}
For almost every $(\alpha,\beta) \in (0,1) \times (0,1)$, the
subshift $\Omega_{\alpha,\beta}$ has infinite index.
\end{prop}

\noindent\textit{Remarks.} (a) This implies that for almost every
$(\alpha,\beta) \in (0,1) \times (0,1)$, $\Omega_{\alpha,\beta}$
is not linearly recurrent.\\[1mm]
(b) Our proof is an extension of arguments from \cite{DP,Kam}.

\begin{proof}
It suffices to show that for each fixed $\beta \in (0,1)$,
$\Omega_{\alpha,\beta}$ has infinite index for almost every
$\alpha \in (0,1)$.

For a sequence $l_k \to \infty$ with
\begin{equation}\label{applkcond}
\sum_{k = 1}^\infty l_k^{-1} = \infty
\end{equation}
(e.g., $l_k = k$), we define the sets $G_{\alpha,\beta}(k)
\subseteq [0,1)$ by
$$
G_{\alpha,\beta}(k) = \{ \theta \in [0,1) : V_\theta (m q_k + j) =
V_\theta(j), \, -2l_k + 1 \le m \le 2l_k - 1, 1 \le j \le q_k \},
$$
where
$$
V_\theta(n) = \chi_{[0,\beta)}(R_\alpha^n \theta).
$$
It is clearly sufficient to show that for each $\beta \in (0,1)$
fixed (and $|\cdot |$ denoting Lebesgue measure),
$$
\left| \limsup_{k \to \infty} G_{\alpha,\beta}(k) \right| > 0
\text{ for almost every $\alpha$}.
$$
Since $| \limsup G_{\alpha,\beta}(k) | \ge \limsup
|G_{\alpha,\beta}(k)|$, this will follow from
\begin{equation}\label{appagoal}
\limsup_{k \rightarrow \infty} |G_{\alpha,\beta}(k)| > 0 \text{
for almost every $\alpha$}.
\end{equation}

Define
\begin{align*}
G^{(1)}_{\alpha,\beta} (k) & = \left\{ \theta : \min_{(-2 l_k + 2)
q_k + 1 \le m \le (2l_k - 1)q_k} \| m \alpha + \theta \| > |q_k
\alpha - p_k|
\right\}, \\
G^{(2)}_{\alpha,\beta} (k) & = \left\{ \theta : \min_{(-2l_k +
2)q_k + 1 \le m \le (2 l_k - 1)q_k} \| m \alpha + \theta - \beta
\|
> |q_k \alpha - p_k| \right\}.
\end{align*}
It follows from \eqref{cfeapp1} that
$$
\|((m \pm q_k) \alpha + \theta) - (m \alpha + \theta) \| = |q_k
\alpha - p_k|.
$$
This in turn implies
\begin{equation}\label{gnrep}
G_{\alpha,\beta} (k) \subseteq G^{(1)}_{\alpha,\beta} (k) \cap
G^{(2)}_{\alpha,\beta} (k)
\end{equation}
On the other hand, we have
\begin{align*}
G^{(1)}_{\alpha,\beta} (k)^c & = \bigcup_{m = (-2l_k+2)q_k +
1}^{(2l_k-1)q_k} \{ \theta :
\| m \alpha + \theta \| \le |q_k \alpha - p_k| \},\\
G^{(2)}_{\alpha,\beta} (k)^c & = \bigcup_{m = (-2l_k + 2)q_k +
1}^{(2l_k-1)q_k} \{ \theta : \| m \alpha + \theta - \beta \| \le
|q_k \alpha - p_k| \},
\end{align*}
which, by \eqref{cfeapp2}, gives for $i = 1,2$,
\begin{equation}\label{gnimeas}
| G^{(i)}_{\alpha,\beta} (k)^c | \le 2 q_k (4l_k - 3) | q_k \alpha
- p_k | \le  (8l_k - 6) \frac{q_k}{q_{k+1}} \le
\frac{8l_k}{a_{k+1}}.
\end{equation}

Combining \eqref{gnrep} and \eqref{gnimeas}, we get
$$
\limsup_{k \rightarrow \infty} |G_{\alpha,\beta}(k)| \ge 1 -
\liminf_{k \rightarrow \infty} \frac{16 l_k}{a_{k+1}}.
$$

By our assumption \eqref{applkcond}, we have that $\liminf_{k
\rightarrow \infty} \frac{16 l_k}{a_{k+1}}$ is less than $1/2$,
say, for almost every $\alpha$ \cite[Theorem~30]{Khin}. This shows
\eqref{appagoal} and hence concludes the proof.
\end{proof}

\section{Some Remarks on Inhomogeneous Diophantine
Approximation}\label{pinnerapp}

Let $\alpha \in (0,1)$ be irrational and let $\gamma \in [0,1)$.  The
\textit{two-sided inhomogeneous approximation constant}
$M(\alpha,\gamma)$ is given by
$$
M(\alpha,\gamma) = \liminf_{|n| \to \infty} |n| \cdot \| n\alpha -
\gamma \|,
$$
where $\| \cdot \|$ denotes the distance from the closest integer.
The number $M(\alpha,\gamma)$ turned out to be important in our
study of the Boshernitzan condition for circle map subshifts
corresponding to partitions of the unit circle into at least three
intervals; compare Section~\ref{Circle}. In this appendix we
sketch a way to compute $M(\alpha,\gamma)$ which was proposed by
Pinner. For background information, we refer the reader to the
excellent texts by Khinchin \cite{Khin} and Rockett and Sz\"usz
\cite{RS}. We shall present results from \cite{Pinner}. Related
work can be found in Cusick et al.\ \cite{CRS} and Komatsu
\cite{Komatsu}.

The \textit{negative continued fraction expansion} of $\alpha$ is
given by
$$
\alpha = \cfrac{1}{a_1 - \cfrac{1}{a_2 - \cfrac{1}{a_3 - \cdots}}}
=: [0;a_1,a_2,a_3,\ldots]^-,
$$
where the integers $a_i \ge 2$ are generated as follows:
$$
\alpha_0 := \{ \alpha \} , \; a_{n+1} := \left\lceil
\frac{1}{\alpha_n} \right\rceil , \; \alpha_{n+1} := \left\lceil
\frac{1}{\alpha_n} \right\rceil - \frac{1}{\alpha_n}.
$$
The corresponding \textit{convergents} $p_n/q_n =
[0;a_1,a_2,\ldots,a_n]^-$ are given by
\begin{alignat*}{3}
p_{-1} &= -1, &\quad        p_0 &= 0,   &\quad      p_{n+1} &= a_{n+1} p_n - p_{n-1},\\
q_{-1} &= 0, &             q_0 &= 1, &           q_{n+1} &=
a_{n+1} q_n - q_{n-1}.
\end{alignat*}
There is a simple way to switch back and forth between regular and
negative continued fraction expansion; see \cite{Pinner}.

Write
$$
\overline{\alpha}_i := [0;a_i,a_{i-1},\ldots,a_1]^- , \;\;
\alpha_i := [0;a_{i+1},a_{i+2},\ldots]^-.
$$
Then
$$
D_i := q_i \alpha - p_i = \alpha_0 \cdots \alpha_i , \; \; q_i =
(\overline{\alpha}_1 \cdots \overline{\alpha}_i)^{-1}.
$$

The $\alpha$-\textit{expansion} of $\gamma$ is now obtained as
follows: Let
$$
\gamma_0 := \{ \gamma \} , \; \; b_{n+1} := \left\lfloor
\frac{\gamma_n}{\alpha_n} \right\rfloor , \; \; \gamma_{n+1} :=
\left\{ \frac{\gamma_n}{\alpha_n} \right\},
$$
so that
$$
\{ \gamma \}  = \sum_{i=1}^\infty b_i D_{i-1}.
$$
Finally, with $t_k := 2b_k - a_k + 2$, let
$$
d_k^- := \sum_{j = 1}^k \frac{t_j \, q_{j-1}}{q_k} , \;\; d_k^+ :=
\sum_{j = k+1}^\infty \frac{t_j \, D_{j-1}}{D_{k-1}}
$$
and
\begin{align*}
s_1 (k) & := \tfrac14 ( 1 - \overline{\alpha}_k + d_k^- ) ( 1 -
\alpha_k + d_k^+ ) q_k D_{k-1}, \\
s_2 (k) & := \tfrac14 ( 1 + \overline{\alpha}_k + d_k^- ) ( 1 +
\alpha_k + d_k^+ ) q_k D_{k-1}, \\
s_3 (k) & := \tfrac14 | 1 - \overline{\alpha}_k - d_k^- | \; | 1
- \alpha_k - d_k^+ | q_k D_{k-1}, \\
s_4 (k) & := \tfrac14 ( 1 + \overline{\alpha}_k - d_k^- ) ( 1 +
\alpha_k + d_k^+ ) q_k D_{k-1}.
\end{align*}

We have the following result \cite{Pinner}:

\begin{theorem}[Pinner]
Suppose that $\gamma \not\in \Z \alpha + \Z$ and that its
$\alpha$-expansion has $b_i = a_i -1$ at most finitely many times.
Then
$$
M(\alpha,\gamma) = \liminf_{ k \to \infty } \min \{ s_1 (k), s_2
(k), s_3 (k), s_4 (k) \}.
$$
\end{theorem}

\end{appendix}

\end{document}